\documentclass[11pt]{article}
\usepackage{amsmath, upgreek}
\usepackage{amssymb}
\usepackage{amscd}
\usepackage{xspace}
\usepackage{verbatim}
\newcommand{\mysection}[1]{
\section{#1}\setcounter{equation}{0}}
\title{\bf  Quasilinear Lane-Emden equations with absorption and measure data}
\author{{\bf Marie-Fran\c{c}oise Bidaut-V\'eron\thanks{E-mail address: veronmf@univ-tours.fr}}\\[1mm]
 {\bf Nguyen Quoc Hung\thanks{ E-mail address: Hung.Nguyen-Quoc@lmpt.univ-tours.fr}}\\[1mm]
 {\bf Laurent V\'eron\thanks{ E-mail address: Laurent.Veron@lmpt.univ-tours.fr}}\\[2mm]
{\small Laboratoire de Math\'ematiques et Physique Th\'eorique, }\\
{\small  Universit\'e Fran\c{c}ois Rabelais,  Tours,  FRANCE}}
\date{}
\begin{document}
 \maketitle


\newcommand{\txt}[1]{\;\text{ #1 }\;}
\newcommand{\tbf}{\textbf}
\newcommand{\tit}{\textit}
\newcommand{\tsc}{\textsc}
\newcommand{\trm}{\textrm}
\newcommand{\mbf}{\mathbf}
\newcommand{\mrm}{\mathrm}
\newcommand{\bsym}{\boldsymbol}
\newcommand{\scs}{\scriptstyle}
\newcommand{\sss}{\scriptscriptstyle}
\newcommand{\txts}{\textstyle}
\newcommand{\dsps}{\displaystyle}
\newcommand{\fnz}{\footnotesize}
\newcommand{\scz}{\scriptsize}
\newcommand{\be}{\begin{equation}}
\newcommand{\bel}[1]{\begin{equation}\label{#1}}
\newcommand{\ee}{\end{equation}}
\newcommand{\eqnl}[2]{\begin{equation}\label{#1}{#2}\end{equation}}
\newcommand{\barr}{\begin{eqnarray}}
\newcommand{\earr}{\end{eqnarray}}
\newcommand{\bars}{\begin{eqnarray*}}
\newcommand{\ears}{\end{eqnarray*}}
\newcommand{\nnu}{\nonumber \\}
\newtheorem{subn}{\name}
\renewcommand{\thesubn}{}
\newcommand{\bsn}[1]{\def\name{#1}\begin{subn}}
\newcommand{\esn}{\end{subn}}
\newtheorem{sub}{\name}[section]
\newcommand{\dn}[1]{\def\name{#1}}   
\newcommand{\bs}{\begin{sub}}
\newcommand{\es}{\end{sub}}
\newcommand{\bsl}[1]{\begin{sub}\label{#1}}
\newcommand{\bth}[1]{\def\name{Theorem}
\begin{sub}\label{t:#1}}
\newcommand{\blemma}[1]{\def\name{Lemma}
\begin{sub}\label{l:#1}}
\newcommand{\bcor}[1]{\def\name{Corollary}
\begin{sub}\label{c:#1}}
\newcommand{\bdef}[1]{\def\name{Definition}
\begin{sub}\label{d:#1}}
\newcommand{\bprop}[1]{\def\name{Proposition}
\begin{sub}\label{p:#1}}
\newcommand{\R}{\eqref}
\newcommand{\rth}[1]{Theorem~\ref{t:#1}}
\newcommand{\rlemma}[1]{Lemma~\ref{l:#1}}
\newcommand{\rcor}[1]{Corollary~\ref{c:#1}}
\newcommand{\rdef}[1]{Definition~\ref{d:#1}}
\newcommand{\rprop}[1]{Proposition~\ref{p:#1}}
\newcommand{\BA}{\begin{array}}
\newcommand{\EA}{\end{array}}
\newcommand{\BAN}{\renewcommand{\arraystretch}{1.2}
\setlength{\arraycolsep}{2pt}\begin{array}}
\newcommand{\BAV}[2]{\renewcommand{\arraystretch}{#1}
\setlength{\arraycolsep}{#2}\begin{array}}
\newcommand{\BSA}{\begin{subarray}}
\newcommand{\ESA}{\end{subarray}}
\newcommand{\BAL}{\begin{aligned}}
\newcommand{\EAL}{\end{aligned}}
\newcommand{\BALG}{\begin{alignat}}
\newcommand{\EALG}{\end{alignat}}
\newcommand{\BALGN}{\begin{alignat*}}
\newcommand{\EALGN}{\end{alignat*}}
\newcommand{\note}[1]{\textit{#1.}\hspace{2mm}}
\newcommand{\Proof}{\note{Proof}}
\newcommand{\qeda}{\hspace{10mm}\hfill $\square$}
\newcommand{\qed}{\\
${}$ \hfill $\square$}
\newcommand{\Remark}{\note{Remark}}
\newcommand{\modin}{$\,$\\[-4mm] \indent}
\newcommand{\forevery}{\quad \forall}
\newcommand{\set}[1]{\{#1\}}
\newcommand{\setdef}[2]{\{\,#1:\,#2\,\}}
\newcommand{\setm}[2]{\{\,#1\mid #2\,\}}
\newcommand{\mt}{\mapsto}
\newcommand{\lra}{\longrightarrow}
\newcommand{\lla}{\longleftarrow}
\newcommand{\llra}{\longleftrightarrow}
\newcommand{\Lra}{\Longrightarrow}
\newcommand{\Lla}{\Longleftarrow}
\newcommand{\Llra}{\Longleftrightarrow}
\newcommand{\warrow}{\rightharpoonup}
\newcommand{
\paran}[1]{\left (#1 \right )}
\newcommand{\sqbr}[1]{\left [#1 \right ]}
\newcommand{\curlybr}[1]{\left \{#1 \right \}}
\newcommand{\abs}[1]{\left |#1\right |}
\newcommand{\norm}[1]{\left \|#1\right \|}
\newcommand{
\paranb}[1]{\big (#1 \big )}
\newcommand{\lsqbrb}[1]{\big [#1 \big ]}
\newcommand{\lcurlybrb}[1]{\big \{#1 \big \}}
\newcommand{\absb}[1]{\big |#1\big |}
\newcommand{\normb}[1]{\big \|#1\big \|}
\newcommand{
\paranB}[1]{\Big (#1 \Big )}
\newcommand{\absB}[1]{\Big |#1\Big |}
\newcommand{\normB}[1]{\Big \|#1\Big \|}
\newcommand{\produal}[1]{\langle #1 \rangle}

\newcommand{\thkl}{\rule[-.5mm]{.3mm}{3mm}}
\newcommand{\thknorm}[1]{\thkl #1 \thkl\,}
\newcommand{\trinorm}[1]{|\!|\!| #1 |\!|\!|\,}
\newcommand{\bang}[1]{\langle #1 \rangle}
\def\angb<#1>{\langle #1 \rangle}
\newcommand{\vstrut}[1]{\rule{0mm}{#1}}
\newcommand{\rec}[1]{\frac{1}{#1}}
\newcommand{\opname}[1]{\mbox{\rm #1}\,}
\newcommand{\supp}{\opname{supp}}
\newcommand{\dist}{\opname{dist}}
\newcommand{\myfrac}[2]{{\displaystyle \frac{#1}{#2} }}
\newcommand{\myint}[2]{{\displaystyle \int_{#1}^{#2}}}
\newcommand{\mysum}[2]{{\displaystyle \sum_{#1}^{#2}}}
\newcommand {\dint}{{\displaystyle \myint\!\!\myint}}
\newcommand{\q}{\quad}
\newcommand{\qq}{\qquad}
\newcommand{\hsp}[1]{\hspace{#1mm}}
\newcommand{\vsp}[1]{\vspace{#1mm}}
\newcommand{\ity}{\infty}
\newcommand{\prt}{\partial}
\newcommand{\sms}{\setminus}
\newcommand{\ems}{\emptyset}
\newcommand{\ti}{\times}
\newcommand{\pr}{^\prime}
\newcommand{\ppr}{^{\prime\prime}}
\newcommand{\tl}{\tilde}
\newcommand{\sbs}{\subset}
\newcommand{\sbeq}{\subseteq}
\newcommand{\nind}{\noindent}
\newcommand{\ind}{\indent}
\newcommand{\ovl}{\overline}
\newcommand{\unl}{\underline}
\newcommand{\nin}{\not\in}
\newcommand{\pfrac}[2]{\genfrac{(}{)}{}{}{#1}{#2}}

\def\ga{\alpha}     \def\gb{\beta}       \def\gg{\gamma}
\def\gc{\chi}       \def\gd{\delta}      \def\ge{\epsilon}
\def\gth{\theta}                         \def\vge{\varepsilon}
\def\gf{\phi}       \def\vgf{\varphi}    \def\gh{\eta}
\def\gi{\iota}      \def\gk{\kappa}      \def\gl{\lambda}
\def\gm{\mu}        \def\gn{\nu}         \def\gp{\pi}
\def\vgp{\varpi}    \def\gr{\rho}        \def\vgr{\varrho}
\def\gs{\sigma}     \def\vgs{\varsigma}  \def\gt{\tau}
\def\gu{\upsilon}   \def\gv{\vartheta}   \def\gw{\omega}
\def\gx{\xi}        \def\gy{\psi}        \def\gz{\zeta}
\def\Gg{\Gamma}     \def\Gd{\Delta}      \def\Gf{\Phi}
\def\Gth{\Theta}
\def\Gl{\Lambda}    \def\Gs{\Sigma}      \def\Gp{\Pi}
\def\Gw{\Omega}     \def\Gx{\Xi}         \def\Gy{\Psi}

\def\CS{{\mathcal S}}   \def\CM{{\mathcal M}}   \def\CN{{\mathcal N}}
\def\CR{{\mathcal R}}   \def\CO{{\mathcal O}}   \def\CP{{\mathcal P}}
\def\CA{{\mathcal A}}   \def\CB{{\mathcal B}}   \def\CC{{\mathcal C}}
\def\CD{{\mathcal D}}   \def\CE{{\mathcal E}}   \def\CF{{\mathcal F}}
\def\CG{{\mathcal G}}   \def\CH{{\mathcal H}}   \def\CI{{\mathcal I}}
\def\CJ{{\mathcal J}}   \def\CK{{\mathcal K}}   \def\CL{{\mathcal L}}
\def\CT{{\mathcal T}}   \def\CU{{\mathcal U}}   \def\CV{{\mathcal V}}
\def\CZ{{\mathcal Z}}   \def\CX{{\mathcal X}}   \def\CY{{\mathcal Y}}
\def\CW{{\mathcal W}} \def\CQ{{\mathcal Q}}
\def\BBA {\mathbb A}   \def\BBb {\mathbb B}    \def\BBC {\mathbb C}
\def\BBD {\mathbb D}   \def\BBE {\mathbb E}    \def\BBF {\mathbb F}
\def\BBG {\mathbb G}   \def\BBH {\mathbb H}    \def\BBI {\mathbb I}
\def\BBJ {\mathbb J}   \def\BBK {\mathbb K}    \def\BBL {\mathbb L}
\def\BBM {\mathbb M}   \def\BBN {\mathbb N}    \def\BBO {\mathbb O}
\def\BBP {\mathbb P}   \def\BBR {\mathbb R}    \def\BBS {\mathbb S}
\def\BBT {\mathbb T}   \def\BBU {\mathbb U}    \def\BBV {\mathbb V}
\def\BBW {\mathbb W}   \def\BBX {\mathbb X}    \def\BBY {\mathbb Y}
\def\BBZ {\mathbb Z}

\def\GTA {\mathfrak A}   \def\GTB {\mathfrak B}    \def\GTC {\mathfrak C}
\def\GTD {\mathfrak D}   \def\GTE {\mathfrak E}    \def\GTF {\mathfrak F}
\def\GTG {\mathfrak G}   \def\GTH {\mathfrak H}    \def\GTI {\mathfrak I}
\def\GTJ {\mathfrak J}   \def\GTK {\mathfrak K}    \def\GTL {\mathfrak L}
\def\GTM {\mathfrak M}   \def\GTN {\mathfrak N}    \def\GTO {\mathfrak O}
\def\GTP {\mathfrak P}   \def\GTR {\mathfrak R}    \def\GTS {\mathfrak S}
\def\GTT {\mathfrak T}   \def\GTU {\mathfrak U}    \def\GTV {\mathfrak V}
\def\GTW {\mathfrak W}   \def\GTX {\mathfrak X}    \def\GTY {\mathfrak Y}
\def\GTZ {\mathfrak Z}   \def\GTQ {\mathfrak Q}

\font\Sym= msam10 
\def\SYM#1{\hbox{\Sym #1}}
\newcommand{\bdw}{\prt\Gw\xspace}
\tableofcontents
\date{}
\maketitle\medskip

\vspace{1mm}
\hspace{.05in}
\medskip

\noindent{\small {\bf Abstract} We study the existence of solutions to the equation $-\Gd_pu+g(x,u)=\gm$
when $g(x,.)$ is a nondecreasing function and $\gm$ a measure. We characterize the good measures, i.e. the ones for which the problem has a renormalized solution. We study particularly the cases where $g(x,u)=\abs x^{-\gb}\abs u^{q-1}u$ and $g(x,u)=\rm{sgn }(u)(e^{\gt\abs u^\gl}-1)$. The results state that a measure is good if it is absolutely continuous with respect to an appropriate Lorentz-Bessel capacities. 
\mysection{Introduction}
Let  $\Gw \sbs \BBR^N$ be a bounded domain containing $0$ and $g:\Gw\ti \BBR \to\BBR$ be a Caratheodory function. We assume that for almost all $x\in\Gw$, $r\mapsto g(x,r)$ is nondecreasing and odd. In this article we consider
the following problem
\bel{E1}\BA {ll}
-\Gd_pu+g(x,u)=\gm\qquad&\text{in }\;\Gw\\
\phantom{-\Gd_p+g(x,u)}
u=0&\text{in }\;\prt\Gw
\EA\ee
where $\Gd_p u=\rm{div}\left(\abs{\nabla u}^{p-2}\nabla u\right)$, ($1<p<N$), is the p-Laplacian and $\gm$ a bounded measure. A measure for which the problem admits a solution, in an appropriate class, is called a {\it good measure}. When $p=2$ and  $g(x,u)=g(u)$ the problem has been considered by Benilan and Brezis \cite {BeBr} in the subcritical case that is when any bounded measure is good. They prove that such is the case if $N\geq 3$ and $g$ satisfies
\bel{E2}\BA {ll}
\myint{1}{\infty}g(s)s^{-\frac{N-1}{N-2}}ds<\infty.
\EA\ee
The supercritical case, always with $p=2$, has been considered by Baras and Pierre \cite{BaPi} when $g(u)=\abs{u}^{q-1}u$ and $q>1$. They prove that the corresponding problem to (\ref{E1}) admits a solution (always unique in that case) if and only if the measure $\gm$ is absolutely continuous with respect to the Bessel capacity $C_{2,q'}$ ($q'=q/(q-1)$). In the case $p\neq 2$ it is shown by Bidaut-V\'eron \cite{Bi2} that if problem (\ref{E1}) with $\gb=0$ and $g(s)=\abs s^{q-1}s$ ($q>p-1>0$) admits a solution, then $\gm$ is absolutely continuous with respect to any capacity $C_{p,\frac{q}{q+1-p}+\ge}$ for any $\ge>0$.\\

In this article we introduce a new class of Bessel capacities which are modelled on Lorentz  spaces $L^{s,q}$ instead of $L^q$ spaces. If $G_\ga$ is the Bessel kernel of order $\ga>0$, we denote by $L^{\ga,s,q}(\BBR^N)$ the Besov space which is the space of functions $\gf=G_\ga\ast f$ for some $f\in L^{s,q}(\BBR^N)$ and we set $\norm\gf_{\ga,s,q}=\norm f_{s,q}$ (a norm which is defined by using rearrangements). Then we set
\bel{E3}\BA {ll}
C_{\ga,s,q}(E)=\inf\{\norm f_{s,q}:\;f\geq 0,\;G_\ga\ast f\geq 1\quad\text{on }E\}
\EA\ee
for any Borel set $E$. We say that a measure $\gm$ in $\Gw$ is absolutely continuous with respect to the capacity $C_{\ga,s,q}$ if , 
\bel{E3}
\forall E\subset \Gw,\,E\text{ Borel },\, C_{\ga,s,q}(E)=0\Longrightarrow \abs\gm(E)=0.
\ee
We also introduce the Wolff potential of a positive measure $\gm\in\mathfrak M_+(\BBR^N)$ by
\bel{E4}\BA {ll}
{\bf W}_{\ga,s}[\gm](x)=\myint{0}{\infty}\left(\myfrac{\gm(B_t(x))}{t^{N-\ga s}}\right)^{\frac{1}{s-1}}\myfrac{dt}{t}
\EA\ee
if $\ga>0$, $1<s<\ga^{-1}N$. When we are dealing with bounded domains $\Gw\subset B_R$ and $\gm\in\mathfrak M_+(\Gw)$, it is useful to introduce truncated Wolff potentials.
\bel{E5}\BA {ll}
{\bf W}^R_{\ga,s}[\gm](x)=\myint{0}{R}\left(\myfrac{\gm(B_t(x))}{t^{N-\ga s}}\right)^{\frac{1}{s-1}}\myfrac{dt}{t}
\EA\ee

We prove the following existence results concerning

\bel{E6}\BA {ll}
-\Gd_pu+\abs x^{-\gb}g(u)=\gm\qquad&\text{in }\;\Gw\\
\phantom{-\Gd_p+\abs x^{-\gb}g(u)}
u=0&\text{in }\;\prt\Gw
\EA\ee

\bth{power} Assume $1<p<N$, $q>p-1$ and  $0\leq\gb <N$ and $\gm$ is a bounded Radon measure in $\Gw$. \smallskip

\noindent 1- If $g(s)=\abs s^{q-1}s$, then (\ref{E6}) admits a renormalized solution if $\gm$ is absolutely continuous with respect to the capacity $C_{p,\frac{Nq}{Nq-(p-1)(N-\gb))},\frac{q}{q+1-p}}$.
\smallskip

\noindent 2- If $g$ satisfies
\bel{}\BA {ll}
\myint{1}{\infty}g(s)s^{-q-1}ds<\infty
\EA\ee
then (\ref{E6}) admits a renormalized solution if $\gm$ is absolutely continuous with respect to the capacity $C_{p,\frac{Nq}{Nq-(p-1)(N-\gb))},1}$.\smallskip

\noindent Furthermore, in both case there holds
\bel{E7}\BA {ll}
-cW_{1,p}^{2\rm{diam\,}(\Gw)}[\gm^-](x)\leq u(x)\leq cW_{1,p}^{2\rm{diam\,}(\Gw)}[\gm^+](x)\qquad\text{for almost all }x\in\Gw.
\EA\ee
where $c$ is a positive constant depending on $p$ and $N$. 
\es

In order to deal with exponential nonlinearities we introduce for $0<\ga<N$ the fractional maximal operator (resp. the truncated fractional maximal operator), defined for a positive measure $\gm$ by
\bel{E8}\BA {ll}\displaystyle
{\bf M}_{\ga}[\gm](x)=\sup_{t>0}\myfrac{\gm(B_t(x))}{t^{N-\ga}},\quad \left(\text{resp }{\bf M}_{\ga,R}[\gm](x)=\sup_{0<t<R}\myfrac{\gm(B_t(x))}{t^{N-\ga}}\right),
\EA\ee
and the $\eta$-fractional maximal operator (resp. the truncated $\eta$-fractional maximal operator)
\bel{E9}\BA {ll}\displaystyle
{\bf M}^\eta_{\ga}[\gm](x)=\sup_{t>0}\myfrac{\gm(B_t(x))}{t^{N-\ga}h_\eta(t)},\quad \left(\text{resp }{\bf M}^\eta_{\ga,R}[\gm](x)=\sup_{0<t<R}\myfrac{\gm(B_t(x))}{t^{N-\ga}h_\eta(t)}\right),
\EA\ee
where $\eta\geq 0$ and
\bel{E10}\BA {ll}
h_\eta(t)=\left\{\BA {ll}(-\ln t)^{-\eta}\qquad&\text{if }0<t<\frac{1}{2}\\[1mm]
(\ln 2)^{-\eta}\qquad&\text{if }t\geq \frac{1}{2}
\EA\right.
\EA\ee

\bth{exp}  Assume $1<p<N$, $\gt>0$ and $\gl\geq 1$. Then there exists $M>0$ depending on $N,p,\gt$ and $\gl$ such that if a measure in $\Gw$, $\gm=\gm^+-\gm^-$ can be decomposed as follows 
\bel{E11}\BA {ll}
\gm^+=f_1+\gn_1\qquad\text{and }\;\gm^-=f_2+\gn_2,
\EA\ee
where  $f_j\in L_+^1(\Gw)$ and $\gn_j\in \mathfrak M_+^b(\Gw)$ ($j=1,2$), and if 
\bel{E12}\BA {ll}
\norm{{\bf M}^{\frac{(p-1)(\gl-1)}{\gl}}_{p,2\rm{diam\,}(\Gw)}[\gn_j]}_{L^\infty(\Gw)}< M,
\EA\ee
there exists a renormalized solution to 
\bel{E13}\BA {ll}
-\Gd_pu+\rm{sign}(u)\left(e^{\gt\abs u^{\gl}}-1\right)=\gm\qquad&\text{in }\;\Gw\\
\phantom{-\Gd_p+\rm{sgn}(s)\left(e^{\gt\abs s^{\gl}}-1\right)}
u=0&\text{in }\;\prt\Gw.
\EA\ee
and satisfies \eqref{E7}. 
\es

Our study is based upon delicate estimates on Wolff potentials and $\eta$-fractional maximal operators which are developed in the first part of this paper.
\section{Lorentz spaces and capacities}
\setcounter{equation}{0}
\subsection{Lorentz spaces}
Let $(X,\Gs,\ga)$ be a measured space. If $f:X\to\BBR$ is a measurable function, we set $S_f(t):=\{x\in X: |f|(x)>t\}$ and $\gl_f(t)=\ga(S_f(t))$. The decreasing rearrangement $f^*$ of $f$ is defined by
$$f^*(t)=\inf\{s>0:\gl_f(s)\leq t\}.
$$
It is well known that $(\Gf(f))^*=\Gf(f^*)$ for any continuous and nondecreasing function $\Gf:\BBR_+\to \BBR_+$. We set
$$f^{**}(t)=\myfrac{1}{t}\myint{0}{t} f^*(\gt)d\gt\qquad\forall t>0.
$$
and, for $1\leq s<\infty$ and $1<q\leq\infty$, 
\bel{L1}
\norm f_{L^{s,q}}=\left\{\BA {ll}\left(\myint{0}{\infty}t^{\frac{q}{s}}(f^{**}(t))^q\myfrac{dt}{t}\right)^{\frac{1}{q}}\qquad&\text{if }q<\infty\\[4mm]
\displaystyle\sup_{\phantom{-}t>0}{\!\!\rm ess}\,t^{\frac{1}{s}}f^{**}(t)\qquad&\text{if }q=\infty
\EA
\right.\ee
It is known that $L^{s,q}(X,\ga)$ is a Banach space when endowed with the norm $\norm ._{L^{s,q}}$. Furthermore
there holds (see e.g. \cite{Gr})
\bel{L2}
\norm{t^{\frac{1}{s}}f^*}_{L^q(\BBR^+,\frac{dt}{t})}\leq \norm f_{L^{s,q}}\leq \myfrac{s}{s-1}\norm{t^{\frac{1}{s}}f^*}_{L^q(\BBR^+,\frac{dt}{t})},
\ee
the left-hand side inequality being valid only if $s>1$. Finally, if $f\in L^{s,q}(\BBR^N)$ (with $1\leq q,s<\infty$ and $\ga$ being the Lebesgue measure) and if $\{\gr_n\}\subset C^{\infty}_c(\BBR^N)$ is a sequence of mollifiers, $f\ast \gr_n\to f$ and $(f\chi_{_{B_n}})\ast \gr_n\to f$ in $L^{s,q}(\BBR^N)$, where $\chi_{_{B_n}}$ is the indicator function of the ball $B_n$ centered at the origin of radius $n$. In particular $C^{\infty}_c(\BBR^N)$ is dense in  $L^{s,q}(\BBR^N)$.
\subsection{Wolff potentials, fractional and $\eta$-fractional maximal operators}

If $D$ is either a bounded domain or whole $\BBR^N$, we denote by $\mathfrak M(D)$ (resp $\mathfrak M^b(D)$) the set of Radon measure (resp. bounded Radon measures) in $D$. Their positive cones are $\mathfrak M_+(D)$ and  $\mathfrak M_+^b(D)$ respectively. If $0<R\leq \infty$ and $\gm\in \mathfrak M_+(D)$ and $R\geq {\rm diam}\,(D)$, we define, for $\ga>0$ and $1<s<\ga^{-1}N$, the $R$-truncated Wolff-potential  by
\bel{L3}
{\bf W}^R_{\ga,s}[\gm](x)=\myint{0}{R}\left(\myfrac{\gm(B_t(x))}{t^{N-\ga s}}\right)^{\frac{1}{s-1}}\myfrac{dt}{t}\quad\text{for a.e. } x\in\BBR^N.
\ee
If $h_\eta(t)=\min\{(-\ln t)^{-\eta},(\ln 2)^{-\eta}\}$ and $0<\ga <N$,  the truncated $\eta$-fractional maximal operator is
\bel{L4}
{\bf M}^\eta_{\ga,R}[\gm](x)=\sup_{0<t<R}\myfrac{\gm(B_t(x))}{t^{N-\ga}h_\eta(t)}\quad\text{for a.e. } x\in\BBR^N.
\ee
If $R=\infty$, we drop it in expressions (\ref{L3}) and (\ref{L4}). In particular
\bel{L5}
\gm(B_t(x))\leq t^{N-\ga}h_\eta(t){\bf M}^\eta_{\ga,R}[\gm](x).
\ee
We also define ${\bf G_\ga}$ the Bessel potential of a measure $\gm$ by
\bel{L6}
{\bf G_\ga}[\gm](x)=\myint{\BBR^N}{}G_\ga(x-y)d\gm(y)\qquad\forall x\in\BBR^N,
\ee
where $G_\ga$ is the Bessel kernel of order $\ga$ in $\BBR^N$.
\bdef{Lorentz} We denote by $L^{\ga,s,q}(\BBR^N)$ the Besov space the space of functions $\gf=G_\ga\ast f$ for some $f\in L^{s,q}(\BBR^N)$ and we set $\norm\gf_{\ga,s,q}=\norm f_{s,q}$. If we set
\bel{E3}\BA {ll}
C_{\ga,s,q}(E)=\inf\{\norm f_{s,q}:\;f\geq 0,\;G_\ga\ast f\geq 1\quad\text{on }E\},
\EA\ee
then $C_{\ga,s,q}$ is a capacity, see \cite{AdHe}.
\es

\subsection{Estimates on potentials}
In the sequel, we denote by $|A|$ the N-dimensional Lebesgue measure of a measurable set $A$ and, if $F,G$ are  functions defined in $\BBR^N$, we set 
$\left\{F>a\right\}:=\{x\in\BBR^N:F(x)>a\}$, $\left\{G\leq b\right\}:=\{x\in\BBR^N:G(x)\leq b\}$ and $\left\{F>a,G\leq b\right\}:=
\left\{F>a\right\}\cap \left\{G\leq b\right\}$. 
The following result is an extension of \cite[Th 1.1]{HoJa}

\bprop{lambda} Let $0\leq\eta<p-1$, $0<\ga p<N$ and $r>0$. There exist $c_0>0$ depending on $N,\ga,p,\eta$ and $\ge_0>0$ depending on $N,\ga,p,\eta, r$ such that, for all $\gm\in\GTM_+(\BBR^N)$ with $diam(supp(\gm))\leq r$ and $R\in (0,\infty]$, $\ge\in (0,\ge_0]$, $\gl>\left(\gm (\BBR^N)\right)^{\frac{1}{p-1}}l(r,R)$ there holds,
\bel{L6}\BA {ll}
\abs{\left\{
{\bf W}^R_{\ga,p}[\gm]>3\gl,({\bf M}^\eta_{\ga p,R}[\gm])^{\frac{1}{p-1}}\leq\ge\gl\right\}}\\[2mm]
\phantom{-----}
\leq c_0\exp\left(-\left(\frac{p-1-\eta}{4(p-1)}\right)^{\frac{p-1}{p-1-\eta}}\ga p\ln 2\,\ge^{-\frac{p-1}{p-1-\eta}}\right)\abs{
\{{\bf W}^R_{\ga,p}[\gm]>\gl\}}.
\EA\ee
where $l(r,R)=\frac{
N-\ga p}{p-1} \left(\min\{r, R\}^{-\frac{N-\ga p}{p-1}}-R^{-\frac{N-\ga p}{p-1}}\right) $ if $R<\infty$, $l(r,R)=\frac{
N-\ga p}{p-1} r^{-\frac{N-\ga p}{p-1}} $ if $R=\infty$. 
Furthermore, if $\eta=0$,  $\ge_0$ is independent of $r$ and \eqref{L6} holds for all $\gm\in\GTM_+(\BBR^N)$ with compact support in $\BBR^N$ and $R\in (0,\infty]$, $\ge\in (0,\ge_0]$, $\gl>0$.
\es
\Proof {\it Case} $R=\infty$. Let $\gl >0$; since ${\bf W}_{\ga,p}[\gm]$ is lower semicontinuous, the set $$D_{\gl}:=\{ {\bf W}_{\ga,p}[\gm]>\gl\}$$ is open. By Whitney covering lemma, there exists a countable set of closed cubes $\{Q_i\}_i$ such that 
$D_{\gl}=\cup_iQ_i$, $\overset{o}{Q_i}\cap \overset{o}{Q_j}=\emptyset$ for $i\neq j$ and
$$\rm {diam} (Q_i)\leq \dist(Q_i,D^c_{\gl})\leq 4\,\rm {diam} (Q_i).
$$
Let $\ge>0$ and $ F_{\ge,\gl}=\left\{
{\bf W}_{\ga,p}[\gm]>3\gl,({\bf M}^\eta_{\ga p}[\gm])^{\frac{1}{p-1}}\leq\ge\gl\right\}$. We claim that there exist $c_0=c_0(N,\ga, p,\eta)>0$ and $\ge_0=\ge_0(N,\ga, p,\eta,r)>0$ such that for any $Q\in \{Q_i\}_i$,  $\ge\in (0,\ge_0]$ and  $\gl>\left(\gm (\BBR^N)\right)^{\frac{1}{p-1}}l(r,\infty)$ there holds
\bel{A0}
\abs{F_{\ge,\gl}\cap Q}\leq c_0 \exp\left(-\left(\frac{p-1-\eta}{4(p-1)}\right)^{\frac{p-1}{p-1-\eta}}\ge^{-\frac{p-1}{p-1-\eta}}\ga p\ln 2\right)\abs Q.
\ee
The first we show that there exists $c_1>0$  depending on $N,\ga, p$ and $\eta$ such that for any $Q\in \{Q_i\}_i$ there holds 
\bel{A3}
F_{\ge,\gl}\cap Q\subset E_{\ge,\gl}~\forall \ge \in (0,c_1], \gl >0
\ee
where 
\bel{A4}
E_{\ge,\gl}=\left\{x\in Q:{\bf W}^{5 \,\rm {diam} (Q)}_{\ga,p}[\gm](x)>\gl,(M_{\ga p}^\eta[\gm](x))^{\frac{1}{p-1}}\leq\ge\gl\right\}.
\ee
Infact, take $Q\in \{Q_i\}_i$ such that $Q\cap F_{\ge,\gl}\neq\emptyset$ and let $x_Q\in D^c_{\gl}$ such that $\dist (x_Q,Q)\leq 4 \,\rm {diam} (Q)$ and ${\bf W}_{\ga,p}[\gm](x_Q)\leq \gl$. For $k\in\BBN$, $r_0=5 \,\rm {diam} (Q)$ and $x\in F_{\ge,\gl}\cap Q$, we have
$$\BA {ll}
\myint{2^kr_0}{2^{k+1}r_0}\left(\myfrac{\gm(B_t(x))}{t^{N-\ga p}}\right)^{\frac {1}{p-1}}\myfrac{dt}{t}
=A+B
\EA$$
where\\
$~~~~~~A=\myint{2^{k}r_0}{2^k\frac{1+2^{k+1}}{1+2^{k}}r_0}\left(\myfrac{\gm(B_t(x))}{t^{N-\ga p}}\right)^{\frac {1}{p-1}}\myfrac{dt}{t}
$
and
$B=\myint{2^k\frac{1+2^{k+1}}{1+2^{k}}r_0}{2^{k+1}r_0}\left(\myfrac{\gm(B_t(x))}{t^{N-\ga p}}\right)^{\frac {1}{p-1}}\myfrac{dt}{t}.
$\\
Since
\bel{A1}\gm(B_t(x))\leq t^{N-\ga p}h_\eta(t)M_{\ga p}^\eta[\gm](x)\leq t^{N-\ga p}h_\eta(t)(\ge\gl)^{p-1}.\ee
Then 
$$B\leq \myint{2^k\frac{1+2^{k+1}}{1+2^{k}}r_0}{2^{k+1}r_0}\left(\myfrac{t^{N-\ga p}h_\eta(t)(\ge\gl)^{p-1}}{t^{N-\ga p}}\right)^{\frac {1}{p-1}}\myfrac{dt}{t}=\ge\gl\myint{2^k\frac{1+2^{k+1}}{1+2^{k}}r_0}{2^{k+1}r_0}\left(h_\eta(t)\right)^{\frac {1}{p-1}}\myfrac{dt}{t}
$$
Replacing $h_\eta(t)$ by its value 
we obtain $B\leq c_2\ge\gl2^{-k}$ after a lengthy computation where $c_2$ depends only on $p$ and $\eta$. Since $\gd:=(\frac{2^k}{2^k+1})^{\frac{N-\ga p}{p-1}}$, then 
$1-\gd\leq c_32^{-k}$ where $c_3$ depends only on $\frac{N-\ga p}{p-1}$, thus
$$\BA {ll}(1-\gd)A
\leq  c_32^{-k} \myint{2^kr_0}{2^{k+1}r_0}\left(\myfrac{\gm(B_t(x))}{t^{N-\ga p}}\right)^{\frac {1}{p-1}}\myfrac{dt}{t}\\[4mm]\phantom{(1-\gd)A}
 \leq c_32^{-k} \ge\gl \myint{2^kr_0}{2^{k+1}r_0}\left(h_\eta(t)\right)^{\frac {1}{p-1}}\myfrac{dt}{t}
 \\[4mm]\phantom{(1-\gd)A}
 \leq c_42^{-k} \ge\gl,
\EA$$
where $c_4=c_4(N,\ga,p,\eta)>0$.\\
By a change of variables and using that for any $x\in F_{\ge,\gl}\cap Q$ and $t\in [r_0(1+2^k),r_0(1+2^{k+1})]$, $B_{\frac{2^kt}{1+2^k}}(x)\subset B_t(x_Q) $, we get
$$\gd A=\int_{r_0(1+2^k)}^{r_0(1+2^{k+1})}\left(\frac{\gm(B_{\frac{2^kt}{1+2^k}})(x)}{t^{N-\ga p}}\right)^{\frac{1}{p-1}}\frac{dt}{t}\leq 
\int_{r_0(1+2^k)}^{r_0(1+2^{k+1})}\left(\frac{\gm(B_t(x_Q))}{t^{N-\ga p}}\right)^{\frac{1}{p-1}}\frac{dt}{t}.
$$
Therefore
$$\myint{2^kr_0}{2^{k+1}r_0}\left(\myfrac{\gm(B_t(x))}{t^{N-\ga p}}\right)^{\frac {1}{p-1}}\myfrac{dt}{t}\leq c_52^{-k} \ge\gl+\int_{r_0(1+2^k)}^{r_0(1+2^{k+1})}\left(\frac{\gm(B_t(x_Q))}{t^{N-\ga p}}\right)^{\frac{1}{p-1}}\frac{dt}{t},
$$
with $c_5=c_5(N,\ga,p,\eta)>0$. This implies
\bel{A2}
\myint{r_0}{\infty}\left(\myfrac{\gm(B_t(x))}{t^{N-\ga p}}\right)^{\frac {1}{p-1}}\myfrac{dt}{t}
\leq 2c_5 \ge\gl+\int_{2r_0}^{\infty}\left(\frac{\gm(B_t(x_Q))}{t^{N-\ga p}}\right)^{\frac{1}{p-1}}\frac{dt}{t}\leq (1+2c_5\ge)\gl,
\ee
 since ${\bf W}_{\ga,p}[\gm](x_Q)\leq \gl$. If $\ge\in (0,c_1]$ with $c_1=(2c_5)^{-1}$ then
$$\myint{r_0}{\infty}\left(\myfrac{\gm(B_t(x)}{t^{N-\ga p}}\right)^{\frac {1}{p-1}}\myfrac{dt}{t}\leq 2\gl
$$
which implies \eqref{A3}.\\
Now, we let $\gl>\left(\gm (\BBR^N)\right)^{\frac{1}{p-1}}l(r,\infty)$. Let $B_1$ be a ball with  radius $r$ such that $supp(\gm)\subset B_1$. We denote $B_2$ by the ball concentric to $B_1$ with radius $2r$. Since $x\notin B_{2}$,   
$${\bf W}_{\ga,p}[\gm](x)=\myint{r}{\infty}\left(\myfrac{\gm(B_t(x))}{t^{N-\ga p}}\right)^{\frac{1}{p-1}}\myfrac{dt}{t}
\leq \left(\gm (\mathbb{R}^N)\right)^{\frac{1}{p-1}}l(r,\infty).
$$
Thus, we obtain $D_{\gl}\subset B_{2}$. In particular, $r_0=5 \,\rm {diam} (Q)\leq 20 r$.\\
Next we set  $m_0=\frac{\max (1,\ln (40 r))}{\ln 2}$, so that $2^{-m}r_0 \leq 2^{-1}$ if $m\geq m_0$. 
Then for any $x\in E_{\ge,\gl}$
$$\BA {l}
\myint{2^{-m}r_0}{r_0}\left(\frac{\gm(B_t(x))}{t^{N-\ga p}}\right)^{\frac{1}{p-1}}\myfrac{dt}{t}
\leq \ge\gl \myint{2^{-m}r_0}{r_0}(h_\eta (t))^{\frac{1}{p-1}}\myfrac{dt}{t}
\\[4mm]\phantom {\myint{2^{-m}r_0}{r_0}\left(\frac{\gm(B_t(x))}{t^{N-\ga p}}\right)^{\frac{1}{p-1}}\myfrac{dt}{t}}
\leq \ge\gl \myint{2^{-m}r_0}{2^{-m_0}r_0}(-\ln t)^{\frac{-\eta}{p-1}}\myfrac{dt}{t}+\ge\gl\myint{2^{-m_0}r_0}{r_0}(\ln 2)^{\frac{-\eta}{p-1}}\myfrac{dt}{t} 
\\[4mm]\phantom {\myint{2^{-m}r_0}{r_0}\left(\frac{\gm(B_t(x))}{t^{N-\ga p}}\right)^{\frac{1}{p-1}}\myfrac{dt}{t}}
\leq m_0\ge\gl+\frac{(p-1)((m-m_0)\ln 2)^{1-\frac{\eta}{p-1}}}{p-1-\eta}\ge\gl.
\EA$$
For the last inequality we have used $a^{1-\frac{\eta}{p-1}}-b^{1-\frac{\eta}{p-1}}\leq (a-b)^{1-\frac{\eta}{p-1}}$ valid for any 
$a\geq b\geq 0$. Therefore, 
\bel{A5}
\myint{2^{-m}r_0}{r_0}\left(\myfrac{\gm(B_t(x))}{t^{N-\ga p}}\right)^{\frac{1}{p-1}}\myfrac{dt}{t}\leq
\myfrac{2(p-1)}{p-1-\eta}m^{1-\frac{\eta}{p-1}}\ge\gl\qquad\forall m\in\BBN, m>m_0^{\frac{p-1}{p-1-\eta}}.
\ee
Set
$$g_i(x)=\myint{2^{-i}r_0}{2^{-i+1}r_0}\left(\myfrac{\gm(B_t(x))}{t^{N-\ga p}}\right)^{\frac{1}{p-1}}\myfrac{dt}{t},
$$
then 
$$\BA {l}
{\bf W}^{r_0}_{\ga,p}[\gm](x)
\leq \myfrac{2(p-1)}{p-1-\eta}m^{1-\frac{\eta}{p-1}}\ge\gl+{\bf W}^{2^{-m}r_0}_{\ga,p}[\gm](x)
\\[2mm]\phantom{{\bf W}^{r_0}_{\ga,p}[\gm](x)}
\leq 
\myfrac{2(p-1)}{p-1-\eta}m^{1-\frac{\eta}{p-1}}\ge\gl+\displaystyle \sum_{i=m+1}^\infty g_i(x)
\EA$$
for all $m> m_0^{\frac{p-1}{p-1-\eta}}$. We deduce that, for $\gb>0$,
\bel{A6}\BA {l}
\abs{E_{\ge,\gl}}\leq\abs{\left\{x\in Q:\displaystyle \sum_{i=m+1}^\infty g_i(x)>\left(1-\myfrac{2(p-1)}{p-1-\eta}m^{1-\frac{\eta}{p-1}}\ge\right)\gl\right\}}
\\[4mm]\phantom{\abs{E_\ge}}
\leq \abs{\left\{x\in Q:\displaystyle \sum_{i=m+1}^\infty g_i(x)>2^{-\gb(i-m-1)}(1-2^{-\gb})\left(1-\myfrac{2(p-1)}{p-1-\eta}m^{1-\frac{\eta}{p-1}}\ge\right)\gl\right\}}
\\[4mm]\phantom{\abs{E_\ge}}
\leq 
\displaystyle \sum_{i=m+1}^\infty
\abs{\left\{x\in Q: g_i(x)>2^{-\gb(i-m-1)}(1-2^{-\gb})\left(1-\myfrac{2(p-1)}{p-1-\eta}m^{1-\frac{\eta}{p-1}}\ge\right)\gl\right\}}.
\EA
\ee
Next we claim that 
\bel{A7}
\abs{\left\{x\in Q: g_i(x)>s\right\}}\leq \frac{c_6(N,\eta)}{s^{p-1}}2^{-i\ga p}\abs{Q}(\ge \gl)^{p-1}.
\ee
To see that, we pick $x_0\in E_{\ge,\gl}$ and we use the Chebyshev's inequality
$$\BA {l}\abs{\left\{x\in Q: g_i(x)>s\right\}}\leq \myfrac{1}{s^{p-1}}\myint{Q}{}\abs{g_i}^{p-1}dx
\\[4mm]\phantom{\abs{\left\{x\in Q: g_i(x)>s\right\}}}
=\myfrac{1}{s^{p-1}}\myint{Q}{}\left(\myint{r_02^{-i}}{r_02^{-i+1}}\left(\myfrac{\gm(B_t(x))}{t^{N-\ga p}}\right)^{\frac{1}{p-1}}\myfrac{dt}{t}\right)^{p-1}dx
\\[4mm]\phantom{\abs{\left\{x\in Q: g_i(x)>s\right\}}}
\leq \myfrac{1}{s^{p-1}}\myint{Q}{}\myfrac{\gm(B_{r_02^{-i+1}}(x))}{(r_02^{-i})^{N-\ga p}}:=A.
\EA$$
Thanks to Fubini's theorem,  the last term $A$ of the above inequality can be rewritten as
$$\BA {l}
A
=\myfrac{1}{s^{p-1}}\myfrac{1}{(r_02^{-i})^{N-\ga p}}\myint{Q}{}\myint{\BBR^N}{} \chi_{B_{r_02^{-i+1}}(x)}(y)d\gm(y) dx
\\[4mm]\phantom{A}
=\myfrac{1}{s^{p-1}}\myfrac{1}{(r_02^{-i})^{N-\ga p}}\myint{Q+B_{r_02^{-i+1}}(0)}{}\myint{Q}{}\chi_{B_{r_02^{-i+1}}(y)}(x)dxd\gm(y) 
\\[4mm]\phantom{A}
\leq  \myfrac{1}{s^{p-1}}\myfrac{1}{(r_02^{-i})^{N-\ga p}}\myint{Q+B_{r_02^{-i+1}}(0)}{}\abs{B_{r_02^{-i+1}}(y)}d\gm(y)
\\[4mm]\phantom{A}
\leq c_7(N)\myfrac{1}{s^{p-1}}2^{-i\ga p}r_0^{\ga p}\gm(Q+B_{r_02^{-i+1}}(0))
\\[4mm]\phantom{A}
\leq c_7(N)\myfrac{1}{s^{p-1}}2^{-i\ga p}r_0^{\ga p}\gm(B_{r_0(1+2^{-i+1})}(x_0)),
\EA$$
since $Q+B_{r_02^{-i+1}}(0)\subset B_{r_0(1+2^{-i+1})}(x_0)$. Using the fact that $\gm(B_t(x_0))\leq (ln~ 2)^{-\eta}t^{N-\ga p}(\ge \gl)^{p-1}$ for all $t>0$ and $r_0=5\,\rm {diam} (Q)$, we obtain
$$A\leq c_8(N,\eta)\myfrac{1}{s^{p-1}}2^{-i\ga p}r_0^{\ga p}(r_0(1+2^{-i+1}))^{N-\ga p}(\ge \gl)^{p-1}\leq c_9(N,\eta)\myfrac{1}{s^{p-1}}2^{-i\ga p}\abs{Q}(\ge \gl)^{p-1},
$$
which is (\ref{A7}). Consequently, (\ref{A6}) can be rewritten as
\bel{A8}\BA {l}
\abs{E_{\ge,\gl}}\leq \displaystyle \sum_{i=m+1}^\infty\myfrac{c_6(N,\eta)}{\left(2^{-\beta(i-m-1)}(1-2^{-\gb})\left(1-\frac{2(p-1)}{p-1-\eta}m^{1-\frac{\eta}{p-1}}\ge\right)\gl\right)^{p-1}}2^{-i\ga p}(\ge \gl)^{p-1}\abs{Q}
\\[6mm]\phantom{\abs{E_\ge}}
\leq c_6(N,\eta)2^{-(m+1)\ga p}\left(\myfrac{\ge}{1-\frac{2(p-1)}{p-1-\eta}m^{1-\frac{\eta}{p-1}}\ge}\right)^{p-1}\abs{Q}\left(1-2^{-\gb}\right)^{-p+1}
\displaystyle\sum_{i=m+1}^\infty 2^{(\gb (p-1)-\ga p)(i-m-1)}.
\EA\ee
If we choose $\gb=\gb(\ga,p)$ so that $\gb (p-1)-\ga p<0$, we obtain 
\bel{A9}
 \abs{E_{\ge,\gl}}\leq c_{10}2^{-m\ga p}\left(\myfrac{\ge}{1-\frac{2(p-1)}{p-1-\eta}m^{1-\frac{\eta}{p-1}}\ge}\right)^{p-1}\abs{Q}\qquad\forall m> m_0^{\frac{p-1}{p-1-\eta}}
\ee
where $c_{10}=c_{10}(N,\ga, p,\eta)>0$. Put $\ge_0=\min \left\{ {\frac{1}{{\frac{{4(p - 1)}}{{p - 1 - \eta }}m_0 + 1}},{c_1}} \right\}$. For any $\ge\in (0,\ge_0]$ we choose $m\in \BBN$ such that
$$
\left(\frac{p-1-\eta}{2(p-1)}\right)^{\frac{p-1}{p-1-\eta}}\left(\frac{1}{\ge}-1\right)^{\frac{p-1}{p-1-\eta}}-1
<m\leq \left(\frac{p-1-\eta}{2(p-1)}\right)^{\frac{p-1}{p-1-\eta}}\left(\frac{1}{\ge}-1\right)^{\frac{p-1}{p-1-\eta}}.$$
Then
$$\left(\myfrac{\ge}{1-\frac{2(p-1)}{p-1-\eta}m^{1-\frac{\eta}{p-1}}\ge}\right)^{p-1}\leq 1
$$
and
$$
2^{-m\ga p}\leq 2^{\ga p-\ga p\left(\frac{p-1-\eta}{2(p-1)}\right)^{\frac{p-1}{p-1-\eta}}\left(\frac{1}{\ge}-1\right)^{\frac{p-1}{p-1-\eta}}}
\leq 2^{\ga p}\exp\left(-\ga p\ln 2\left(\frac{p-1-\eta}{4(p-1)}\right)^{\frac{p-1}{p-1-\eta}}\ge^{-\frac{p-1}{p-1-\eta}}\right).
$$
Combining these inequalities  with \eqref{A9} and \eqref{A3}, we get (\ref{A0}).

In the case $\eta =0$ we still have for any $m\in\mathbb{N}$, $\gl ,\ge>0$ and  $x\in E_{\ge,\gl}$
$${\bf W}^{r_0}_{\ga,p}[\gm](x)\leq 
m\ge\gl+\displaystyle \sum_{i=m+1}^\infty g_i(x)
$$
Accordingly (\ref{A9}) reads as 
$$\abs{ E_{\ge,\gl}}\leq 
c_{10}2^{-m\ga p}\left(\myfrac{\ge}{1-m\ge}\right)^{p-1}\abs{Q}~~\forall m\in\mathbb{N}, \gl ,\ge>0 \textrm{ with } m\ge <1.
$$
Put $\ge_0=\min\{\frac{1}{2},c_1\}$. For any $\ge\in (0,\ge_0]$ and $m\in\BBN$ satisfies $\ge^{-1}-2<m\leq \ge^{-1}-1$, we finally get from \eqref{A3}
\bel{A12}
\abs{F_{\ge,\gl}\cap Q}\leq \abs{ E_{\ge,\gl}}\leq c_{10}2^{2\ga p}\exp\left(-\ga p\ge^{-1}\ln 2\right)|Q|,
\ee
which ends the proof in the case $R=\infty$.\medskip

\noindent {\it Case} $R<\infty$. For $\gl>0$, $D_{\gl}=\{{\bf W}_{\ga,p}^R>\gl\}$ is open. Using again Whitney covering lemma, there exists a countable set of closed cubes $\CQ:=\{Q_i\}$ such that $\cup_iQ_i=D_{\gl}$, $\overset{o}{Q_i}\cap \overset{o}{Q_j}=\emptyset$ for $i\neq j$ and $ \dist(Q_i,D^c_\gl)\leq 4\,\rm {diam} (Q_i)$. If $Q\in \CQ:$ is such that $\rm{diam}\, (Q)>\frac{R}{8}$, there exists a finite number $n_Q$ of closed dyadic cubes $\{P_{j,Q}\}_{j=1}^{n_Q}$ such that $\cup_{j=1}^{n_Q} P_{j,Q}=Q$, 
$\overset{o}{P_{i,Q}}\cap \overset{o}{P_{j,Q}}=\emptyset$ if $i\neq j$ and 
$\frac{R}{16}<\rm{diam}\, (P_{j,Q})\leq \frac{R}{8}$. We set
$\CQ'=\left\{Q\in\CQ: \rm{diam}\, (Q)\leq\frac{R}{8}\right\}$, $\CQ''= \left\{P_{i,Q}: 1\leq i\leq n_Q,Q\in\CQ,\rm{diam}\, (Q)>\frac{R}{8}\right\}$ and $\CF=\CQ'\cup\CQ''$. \\
For $\ge>0$ we denote again
$F_{\ge,\gl}=\left\{
{\bf W}^R_{\ga,p}[\gm]>3\gl,({\bf M}^\eta_{\ga p,R}[\gm])^{\frac{1}{p-1}}\leq\ge\gl\right\}$. Let $Q\in\CF$ such that $F_{\ge,\gl}\cap Q\neq\emptyset$ and $r_0=5\,\rm{diam}\, (Q)$. \smallskip

\noindent If $\dist(D^c_\gl,Q)\leq 4\,\rm{diam}\, (Q)$, that is if there exists $x_Q\in D^c_\gl$ such that 
$\dist(x_Q,Q)\leq 4\,\rm{diam}\, (Q)$ and ${\bf W}_{\ga,p}^R[\gm](x_Q)\leq \gl$, we find, by the same argument as in the case $R=\infty$, (\ref{A2}), that for any $x\in F_{\ge,\gl}\cap Q$ there holds
\bel{A13}
\myint{r_0}{R}\left(\myfrac{\gm(B_t(x))}{t^{N-\ga p}}\right)^{\frac {1}{p-1}}\myfrac{dt}{t}
\leq (1+c_{11}\ge)\gl.
\ee
where $c_{11}=c_{11}(N,\ga,p,\eta)>0$.\\
\noindent If $\dist(D^c_\gl,Q)> 4\,\rm{diam}\, (Q)$, we have $\frac{R}{16}<\rm{diam}\, (Q)\leq \frac{R}{8}$
 since $Q\in  \CQ''$. Then, for all $x\in F_{\ge,\gl}\cap Q$, there holds
 \bel{A14}\BA {l}
\myint{r_0}{R}\left(\myfrac{\gm(B_t(x))}{t^{N-\ga p}}\right)^{\frac {1}{p-1}}\myfrac{dt}{t} \leq \myint{\frac {5R}{16}}{R}\left(\myfrac{t^{N-\ga p}(\ln 2)^{-\eta}(\ge \gl)^{p-1}}{t^{N-\ga p}}\right)^{\frac {1}{p-1}}\myfrac{dt}{t}
  \\[4mm]\phantom{ \myint{r_0}{R}\left(\myfrac{\gm(B_t(x))}{t^{N-\ga p}}\right)^{\frac {1}{p-1}}\myfrac{dt}{t}}
  =(\ln 2)^{-\frac{\eta}{p-1}}\ln \frac{16}{5}~ \ge \gl
    \\[4mm]\phantom{ \myint{r_0}{R}\left(\myfrac{\gm(B_t(x))}{t^{N-\ga p}}\right)^{\frac {1}{p-1}}\myfrac{dt}{t}}
   \leq 2\ge\gl.
 \EA\ee
Thus, if we take $\ge\in (0,c_{12}]$ with $c_{12}=\min\{1,c_{11}^{-1}\}$, we derive
\bel{A15}F_{\ge,\gl}\cap Q\subset E_{\ge,\gl},
\ee
where $$E_{\ge,\gl}=\left\{ {\bf W}_{\ga,p}^{r_0}[\gm]>\gl, \left({\bf M}_{\ga p,R}^\eta[\gm]\right)^{\frac{1}{p-1}}\leq\ge\gl\right\}.$$
Furthermore, since $x\notin B_{2}$,   
$${\bf W}_{\ga,p}^R[\gm](x)=\myint{\min\{r,R\}}{R}\left(\myfrac{\gm(B_t(x))}{t^{N-\ga p}}\right)^{\frac{1}{p-1}}\myfrac{dt}{t}
\leq \left(\gm (\mathbb{R}^N)\right)^{\frac{1}{p-1}}l(r,R).
$$
Thus, if $\gl>\left(\gm (\BBR^N)\right)^{\frac{1}{p-1}}l(r,R)$ then $D_{\gl}\subset B_{2}$ which implies $r_0=5 \,\rm {diam} (Q)\leq 20 r$.\\
The end of the proof is as in the case $R=\infty$.\qeda\medskip

In the next result we list a series of equivalent norms concerning Radon measures.

\bth{equiv} Assume $\ga>0$, $0<p-1<q<\infty$, $0<\ga p<N$ and $0<s\leq \infty$. Then there exists a constant $c_{13}=c_{13}(N,\ga, p,q,s)>0$ such that for any $R\in (0,\infty]$ and $\gm\in\mathfrak M_+(\BBR^N)$, there holds
\bel{B1}
c_{13}^{-1}\norm{{\bf W}_{\ga,p}^R[\gm]}_{L^{q,s}(\BBR^N)}\leq 
\norm{{\bf M}_{\ga p,R}[\gm]}_{L^{\frac{q}{p-1},\frac{s}{p-1}}(\BBR^N)}^{\frac{1}{p-1}}\leq 
c_{13}\norm{{\bf W}_{\ga,p}^R[\gm]}_{L^{q,s}(\BBR^N)}.
\ee
For any $R>0$, there exists $c_{14}=c_{14}(N,\ga, p,q,s,R)>0$ such that for any $\gm\in\mathfrak M_+(\BBR^N)$,
\bel{B2}
c_{14}^{-1}\norm{{\bf W}_{\ga,p}^R[\gm]}_{L^{q,s}(\BBR^N)}\leq 
\norm{{\bf G}_{\ga p}[\gm]}_{L^{\frac{q}{p-1},\frac{s}{p-1}}(\BBR^N)}^{\frac{1}{p-1}}\leq 
c_{14}\norm{{\bf W}_{\ga,p}^R[\gm]}_{L^{q,s}(\BBR^N)}.
\ee
In (\ref {B2}), $\norm{{\bf W}_{\ga,p}^R[\gm]}_{L^{q,s}(\BBR^N)}$ can be replaced by 
$\norm{{\bf M}_{\ga p,R}[\gm]}_{L^{\frac{q}{p-1},\frac{s}{p-1}}(\BBR^N)}^{\frac{1}{p-1}}$.
\es
\Proof We denote $\gm_n$ by ${\bf \chi}_{B_n}\gm$ for  $n\in\BBN^\ast$.\\
\noindent{\it Step 1} We claim that 
\bel{B3}
\norm{{\bf W}_{\ga,p}^R[\gm]}_{L^{q,s}(\BBR^N)}\leq 
c{'}\!_{13}\norm{{\bf M}_{\ga p,R}[\gm]}_{L^{\frac{q}{p-1},\frac{s}{p-1}}(\BBR^N)}^{\frac{1}{p-1}}.
\ee
From \rprop{lambda} there exist positive constants $c_0=c_0(N,\ga,p),a=a(\ga,p)$ and $\ge_0=\ge_0(N,\ga,p)$ such that for all $n\in\BBN^\ast$, $t>0$, $0<R\leq\infty$ and $0<\ge\leq\ge_0$, there holds
\bel{B4}\BA {ll}
\abs{\left\{
{\bf W}^R_{\ga,p}[\gm_n]>3t,({\bf M}^\eta_{\ga p,R}[\gm_n])^{\frac{1}{p-1}}\leq\ge t\right\}}
\leq c_0\exp\left(-a\ge^{-1}\right)\abs{
\{{\bf W}^R_{\ga,p}[\gm_n]>t\}}.
\EA\ee
In the case $0<s<\infty$ and $0<q<\infty$, we have
$$
\abs{\left\{
{\bf W}^R_{\ga,p}[\gm_n]>3t\right\}}^{\frac{s}{q}}
\leq c_{15}\exp\left(-\frac{s}{q}a\ge^{-1}\right)\abs{
\{{\bf W}^R_{\ga,p}[\gm_n]>t\}}^{\frac{s}{q}}+c_{15}
\abs{\left\{
({\bf M}^\eta_{\ga p,R}[\gm_n])^{\frac{1}{p-1}}>\ge t\right\}}^{\frac{s}{q}}.
$$
with $c_{15}=c_{15}(N,\ga,p,q,s)>0$.\\
Multiplying by $t^{s-1}$ and integrating over $(0,\infty)$, we obtain
$$\BA {l}
\myint{0}{\infty}t^s\abs{\left\{
{\bf W}^R_{\ga,p}[\gm_n]>3t\right\}}^{\frac{s}{q}}\myfrac{dt}{t}\leq
c_{15}\exp\left(-\frac{s}{q}a\ge^{-1}\right)\myint{0}{\infty}t^s
\abs{
\{{\bf W}^R_{\ga,p}[\gm_n]>t\}}^{\frac{s}{q}}\myfrac{dt}{t}
\\[4mm]\phantom{--\myint{0}{\infty}t^s\abs{\left\{
{\bf W}^R_{\ga,p}[\gm_n]>3t\right\}}^{\frac{s}{q}}\myfrac{dt}{t}}
+c_{15}\myint{0}{\infty}t^s\abs{\left\{
{\bf M}^\eta_{\ga p,R}[\gm_n]>(\ge t)^{{p-1}}\right\}}^{\frac{s}{q}}
\myfrac{dt}{t}.
\EA$$
By a change of variable, we derive
$$\BA {l}\left(3^{-s}-c_{15} \exp\left(-\frac{s}{q}a\ge^{-1}\right)\right)\myint{0}{\infty}t^s
\abs{
\{{\bf W}^R_{\ga,p}[\gm_n]>t\}}^{\frac{s}{q}}\myfrac{dt}{t}
\\[4mm]\phantom{\left(3^{-s} \exp\left(-a\ge^{-1}\right)\right)\myint{0}{\infty}t^s
\abs{
\{{\bf W}^R_{\ga,p}>t\}}}
\leq \myfrac{c_{15}\ge^{-s}}{p-1}\myint{0}{\infty}t^{\frac{s}{p-1}}\abs{\left\{
{\bf M}^\eta_{\ga p,R}[\gm_n]>t\right\}}^{\frac{s}{q}}
\myfrac{dt}{t}.
\EA$$
We choose $\ge$ small enough so that $3^{-s}-c_{15} \exp\left(-\frac{s}{q}a\ge^{-1}\right)>0$, we derive from \eqref{L2} and ${\left\| {{t^{1/{s_1}}}{f^*}} \right\|_{{L^{s_2}}\left( {\BBR,\frac{{dt}}{t}} \right)}} = {s_1^{1/{s_2}}}{\left\| {\lambda _f^{1/{s_1}}t} \right\|_{{L^{s_2}}\left( {\BBR,\frac{{dt}}{t}} \right)}}$ for any $f\in L^{s_1,s_2}(\BBR^N)$ with $0<s_1<\infty,0<s_2\leq \infty $ 
$$\norm{{\bf W}^R_{\ga,p}[\gm_n]}_{L^{q,s}(\BBR^N)}\leq c'\!_{13}
\norm{{\bf M}_{\ga p,R}[\gm_n]}_{L^{\frac{q}{p-1},\frac{s}{p-1}}(\BBR^N)}^{\frac{1}{p-1}},
$$
and (\ref{B3}) follows by Fatou's lemma. Similarly, we can prove (\ref{B3}) in the case $s=\infty$.\smallskip

\noindent{\it Step 2} We claim that 
\bel{B6}
\norm{{\bf W}_{\ga,p}^R[\gm]}_{L^{q,s}(\BBR^N)}\geq 
c{''}\!_{13}\norm{{\bf M}_{\ga p,R}[\gm]}_{L^{\frac{q}{p-1},\frac{s}{p-1}}(\BBR^N)}^{\frac{1}{p-1}}.
\ee
For $R>0$ we have 
\bel{B7}\BA {ll}
{\bf W}_{\ga,p}^{2R}[\gm_n](x)={\bf W}_{\ga,p}^{R}[\gm_n](x)
+\myint{R}{2R}\left(\myfrac{\gm_n(B_t(x))}{t^{N-\ga p}}\right)^{\frac{1}{p-1}}\myfrac{dt}{t}
\\[4mm]\phantom{{\bf W}_{\ga,p}^{2R}[\gm](x)}
\leq  {\bf W}_{\ga,p}^{R}[\gm_n](x)+\left(\myfrac{\gm_n(B_{2R}(x))}{R^{N-\ga p}}\right)^{\frac{1}{p-1}}.
\EA\ee
Thus
$$\abs{\left\{x:{\bf W}_{\ga,p}^{2R}[\gm_n](x)>2t\right\}}
\leq \abs{\left\{x:{\bf W}_{\ga,p}^{R}[\gm_n](x)>t\right\}}+\abs{\left\{x:\myfrac{\gm_n(B_{2R}(x))}{R^{N-\ga p}}>t^{p-1}\right\}},
$$
 Consider $\{z_j\}_{i=1}^m\subset B_2$ such that $B_2\subset\bigcup_{i=1}^mB_{\frac{1}{2}}(z_i)$. Thus
$B_{2R}(x)\subset\bigcup_{i=1}^mB_{\frac{R}{2}}(x+Rz_i)$ for any $x\in\BBR^N$ and $R>0$. Then
$$\BA {l}
\abs{\left\{x:\myfrac{\gm_n(B_{2R}(x))}{R^{N-\ga p}}>t^{p-1}\right\}}
\leq \abs{\left\{x:\displaystyle\sum_{i=1}^m\myfrac{\gm_n(B_{\frac{R}{2}}(x+Rz_i))}{R^{N-\ga p}}>t^{p-1}\right\}}
\\[4mm]\phantom{\abs{\left\{x:\myfrac{\gm_n(B_{2R}(x))}{R^{N-\ga p}}>t^{p-1}\right\}}}
\leq \displaystyle\sum_{i=1}^m
\abs{\left\{x:\myfrac{\gm_n(B_{\frac{R}{2}}(x+Rz_i))}{R^{N-\ga p}}>\myfrac{1}{m}t^{p-1}\right\}}
\\[4mm]\phantom{\abs{\left\{x:\myfrac{\gm_n(B_{2R}(x))}{R^{N-\ga p}}>t^{p-1}\right\}}}
\leq \displaystyle\sum_{i=1}^m\abs{\left\{x-Rz_i:\myfrac{\gm_n(B_{\frac{R}{2}}(x))}{R^{N-\ga p}}>\myfrac{1}{m}t^{p-1}\right\}}
\\[4mm]\phantom{\abs{\left\{x:\myfrac{\gm_n(B_{2R}(x))}{R^{N-\ga p}}>c_9t^{p-1}\right\}}}
=m\abs{\left\{x:\myfrac{\gm_n(B_{\frac{R}{2}}(x))}{R^{N-\ga p}}>\myfrac{1}{m}t^{p-1}\right\}}.
\EA$$
Moreover from (\ref{B7}) 
$$\left(\myfrac{\gm_n(B_{\frac{R}{2}} (x))}{R^{N-\ga p}}\right)^{\frac{1}{p-1}}
\leq 2{\bf W}_{\ga,p}^{R}[\gm_n](x),
$$
thus
$$\abs{\left\{x:\myfrac{\gm_n(B_{2R}(x))}{R^{N-\ga p}}>t^{p-1}\right\}}
\leq m\abs{\left\{x:{\bf W}_{\ga,p}^{R}[\gm_n](x)>\frac{1}{2m^{\frac{1}{p-1}}}t\right\}}.
$$
This leads to
$$
\abs{\left\{x:{\bf W}_{\ga,p}^{2R}[\gm_n](x)>2t\right\}}
\leq (m+1)\abs{\left\{x:{\bf W}_{\ga,p}^{R}[\gm_n](x)>\frac{1}{2m^{\frac{1}{p-1}}}t\right\}}~\forall t>0
$$
 This implies
$$
\norm{{\bf W}_{\ga,p}^{2R}[\gm_n]}_{L^{\frac{q}{p-1},\frac{s}{p-1}}(\BBR^N)}
\leq c_{16}\norm{{\bf W}_{\ga,p}^{R}[\gm_n]}_{L^{\frac{q}{p-1},\frac{s}{p-1}}(\BBR^N)}.
$$
with $c_{16}=c_{16}(N,\ga, p,q,s)>0$.
By  Fatou's lemma, we get \bel{B9}\BA {ll}
\norm{{\bf W}_{\ga,p}^{2R}[\gm]}_{L^{\frac{q}{p-1},\frac{s}{p-1}}(\BBR^N)}
\leq c_{16}\norm{{\bf W}_{\ga,p}^{R}[\gm]}_{L^{\frac{q}{p-1},\frac{s}{p-1}}(\BBR^N)}.
\EA\ee
On the other hand, from the identity in (\ref{B7}) we derive that for any $\gr\in (0,R)$,
 $$
{\bf W}_{\ga,p}^{2R}[\gm](x)\geq {\bf W}_{\ga,p}^{2\gr}[\gm](x)\geq c_{17}\sup_{0<\gr\leq R}\left(\myfrac{\gm(B_\gr(x))}{\gr^{N-\ga p}}\right)^{\frac{1}{p-1}},
$$
with $c_{17}=c_{17}(N,\ga, p)>0$, from which follows
 \bel{B11}\BA {ll}
{\bf W}_{\ga,p}^{2R}[\gm](x)\geq c_{17}\left({\bf M}_{\ga p,R}[\gm](x)\right)^{\frac{1}{p-1}}.
\EA\ee
Combining (\ref{B9}) and (\ref{B11}) we obtain (\ref{B6}) and then (\ref{B1}). Notice that the estimates are independent of $R$ and thus valid if $R=\infty$.\smallskip

\noindent{\it Step 3} We claim that (\ref{B2}) holds. By the previous result we have also
\bel{B12}
c_{18}^{-1}\norm{{\bf W}_{\frac{\ga p}{2},2}^R[\gm]}_{L^{\frac{q}{p-1},\frac{s}{p-1}}(\BBR^N)}\leq 
\norm{{\bf M}_{\ga p,R}[\gm]}_{L^{\frac{q}{p-1},\frac{s}{p-1}}(\BBR^N)}\leq 
c_{18}\norm{{\bf W}_{\frac{\ga p}{2},2}^R[\gm]}_{L^{\frac{q}{p-1},\frac{s}{p-1}}(\BBR^N)}.
\ee
where  $c_{18}=c_{18}(N,\ga,p,q,s)>0$.
For $R>0$, the Bessel kernel satisfies\cite[V-3-1]{St}
$$c_{19}^{-1}\left(\myfrac{\chi_{B_R}(x)}{\abs x^{N-\ga p}}\right)\leq G_{\ga p}(x)\leq 
c_{19}\left(\myfrac{\chi_{B_{\frac{R}{2}}}(x)}{\abs x^{N-\ga p}}\right)+c_{19}e^{-\frac{\abs x}{2}}
\qquad\forall x\in\BBR^N,$$
where $c_{19}=c_{19}(N,\ga, p, R)>0$. Therefore
\bel{B12'}c_{19}^{-1}\left(\myfrac{\chi_{B_R}}{\abs .^{N-\ga p}}\right)\ast\gm\leq {\bf G}_{\ga p}[\gm]\leq 
c_{19}\left(\myfrac{\chi_{B_{\frac{R}{2}}}}{\abs .^{N-\ga p}}\right)\ast\gm+c_{19}e^{-\frac{\abs .}{2}}
\ast\gm.
\ee
By integration by parts, we get
$$\left(\myfrac{\chi_{B_R}}{\abs .^{N-\ga p}}\right)\ast\gm (x)=(N-\ga p){\bf W}_{\frac{\ga p}{2},2}^R[\gm](x)+\myfrac{\gm(B_R(x))}{R^{N-\ga p}}\geq (N-\ga p){\bf W}_{\frac{\ga p}{2},2}^R[\gm](x),
$$
which implies
\bel{B13}
c_{20}\norm{{\bf W}_{\frac{\ga p}{2},2}^R[\gm]}_{L^{\frac{q}{p-1},\frac{s}{p-1}}(\BBR^N)}
\leq \norm{{\bf G}_{\ga p}[\gm]}_{L^{\frac{q}{p-1},\frac{s}{p-1}}(\BBR^N)}.
\ee
where  $c_{20}=c_{20}(N,\ga,p,q,s)>0$.
Furthermore $e^{-\frac{\abs x}{2}}\leq c_{21}\chi_{B_{\frac{R}{2}}}\ast e^{-\frac{\abs .}{2}}(x)$ where 
$c_{21}=c_{21}(N,R)>0$, thus
$$e^{-\frac{\abs .}{2}}\ast \gm\leq c_{21}\left(\chi_{B_{\frac{R}{2}}}\ast e^{-\frac{\abs .}{2}}\right)\ast \gm
=c_{21}e^{-\frac{\abs .}{2}}\ast \left(\chi_{B_{\frac{R}{2}}}\ast \gm\right).
$$
Since 
$$\chi_{B_{\frac{R}{2}}}\ast \gm(x)=\gm(B_{\frac{R}{2}}(x))\leq c_{22}{\bf W}_{\frac{\ga p}{2},2}^R[\gm](x)$$
where $c_{22}=c_{22}(N,\ga, p,R)>0$, we derive with $c_{23}=c_{21}c_{22}$
$$e^{-\frac{\abs .}{2}}\ast \gm\leq c_{23}e^{-\frac{\abs .}{2}}\ast{\bf W}_{\frac{\ga p}{2},2}^R[\gm].
$$
Using  Young inequality, we obtain
\bel{B14}\BA {l}
\norm{e^{-\frac{\abs .}{2}}\ast \gm}_{L^{\frac{q}{p-1},\frac{s}{p-1}}(\BBR^N)}
\leq c_{23}\norm{e^{-\frac{\abs .}{2}}\ast{\bf W}_{\frac{\ga p}{2},2}^R[\gm]}_{L^{\frac{q}{p-1},\frac{s}{p-1}}(\BBR^N)}
\\[4mm]\phantom{\norm{e^{-\frac{\abs .}{2}}\ast \gm}_{L^{\frac{q}{p-1},\frac{s}{p-1}}(\BBR^N)}}
\leq c_{24}\norm{{\bf W}_{\frac{\ga p}{2},2}^R[\gm]}_{L^{\frac{q}{p-1},\frac{s}{p-1}}(\BBR^N)}
\norm{e^{-\frac{\abs .}{2}}}_{L^{1,\infty}(\BBR^N)}
\\[4mm]\phantom{\norm{e^{-\frac{\abs .}{2}}\ast \gm}_{L^{\frac{q}{p-1},\frac{s}{p-1}}(\BBR^N)}}
\leq c_{25}\norm{{\bf W}_{\frac{\ga p}{2},2}^R[\gm]}_{L^{\frac{q}{p-1},\frac{s}{p-1}}(\BBR^N)}.
\EA\ee
where $c_{25}=c_{25}(N,\ga, p,R)>0$.\\
Since by integration by parts there holds as above
$$\left(\myfrac{\chi_{B_{\frac{R}{2}}}}{\abs .^{N-\ga p}}\right)\ast\gm (x)=(N-\ga p){\bf W}_{\frac{\ga p}{2},2}^{\frac{R}{2}}[\gm](x)+2^{N-\ga p}\myfrac{\gm(B_{\frac{R}{2}}(x))}{R^{N-\ga p}}\leq c_{26}{\bf W}_{\frac{\ga p}{2},2}^R[\gm](x),
$$
where $c_{26}=c_{26}(N,\ga,p)>0$
we obtain
\bel{B15}
\norm{\left(\myfrac{\chi_{B_R}}{\abs .^{N-\ga p}}\right)\ast\gm}_{L^{\frac{q}{p-1},\frac{s}{p-1}}(\BBR^N)}\leq c_{27}\norm{{\bf W}_{\frac{\ga p}{2},2}^R[\gm]}_{L^{\frac{q}{p-1},\frac{s}{p-1}}}.
\ee
where $c_{27}=c_{27}(N,\ga, p,q,s)>0$.
Thus 
\bel{B16}
\norm{{\bf G}_{\ga p}[\gm]}_{L^{\frac{q}{p-1},\frac{s}{p-1}}(\BBR^N)}\leq c_{28}\norm{{\bf W}_{\frac{\ga p}{2},2}^R[\gm]}_{L^{\frac{q}{p-1},\frac{s}{p-1}}}.
\ee
where $c_{28}=c_{28}(N,\ga, p,q,s,R)>0$.\\
follows by combining \eqref{B12'}, \eqref{B14} and \eqref{B15}. Then, combining \eqref{B13}, \eqref{B16} and using \eqref{B12}, \eqref{B1} we obtain \eqref{B2}.\qeda\medskip

\noindent{\Remark} Proposition 5.1 in \cite{PhVer} is a particular case of the previous result.

\bth{exp1} Let $\ga>0$, $p>1$, $0\leq\eta<p-1$, $0<\ga p<N$ and $r>0$. Set $\gd_0=\left(\frac{p-1-\eta}{12(p-1)}\right)^{\frac{p-1}{p-1-\eta}}\ga p\ln 2$. Then there exists $c_{29}>0$, depending on 
 $N$, $\ga$, $p$, $\eta$  and $r$ such that for any $R\in (0,\infty]$, $\gd\in (0,\gd_0)$, $\gm\in\mathfrak M_+(\BBR^N)$,  any  ball $B_1\subset \BBR ^N$ with radius $\leq r$ and ball $B_2$ concentric to $B_1$ with radius double $B_1$'s radius, there holds
 \bel{B17}\BA {l}
\frac{1}{\abs{B_{2}}}\myint{B_{2}}{}
\exp\left(\gd \frac{\left({\bf W}^R_{\ga,p}[\gm_{B_{1}}](x)\right)^{\frac{p-1}{p-1-\eta}}}{\norm{{\bf M}^\eta_{\ga p,R}[\gm_{B_{1}}]}_{L^\infty(B_{1})}^{\frac{1}{p-1-\eta}}}\right) dx
\leq \myfrac{c_{29}}{\gd_0-\gd}
\EA\ee
where  $\gm_{B_{1}}=\chi _{B_{1}}\gm$. Furthermore, if $\eta=0$, $c_{29}$ is independent of $r$.\es
\Proof Let $\gm\in\mathfrak M_+(\BBR^N)$ such that $M:=\norm{{\bf M}^\eta_{\ga p,R}[\gm_{B_{1}}]}_{L^\infty(B_{1})}<\infty$. By 
\rprop{lambda}-(\ref{L6}) with $\gm=\gm_{B_{1}}$, there exist $c_0>0$ depending on $N,\ga,p,\eta$ and $\ge_0>0$ depending on $N,\ga,p,\eta$ and $ r$ such that, for all  $R\in (0,\infty]$, $\ge\in (0,\ge_0]$, $t>\left(\gm_{B_{1}} (\BBR^N)\right)^{\frac{1}{p-1}}l(r',R)$ where $r'$ is radius  of $B_1$ there holds,
\bel{}\BA {ll}
\abs{\left\{
{\bf W}^R_{\ga,p}[\gm_{B_{1}}]>3t,({\bf M}^\eta_{\ga p,R}[\gm_{B_{1}}])^{\frac{1}{p-1}}\leq\ge t\right\}}\\[2mm]
\phantom{-----}
\leq c_0\exp\left(-\left(\frac{p-1-\eta}{4(p-1)}\right)^{\frac{p-1}{p-1-\eta}}\ga p\ln 2\,\ge^{-\frac{p-1}{p-1-\eta}}\right)\abs{
\{{\bf W}^R_{\ga,p}[\gm_{B_{1}}]>t\}}.
\EA\ee
Since $\left(\gm_{B_{1}} (\BBR^N)\right)^{\frac{1}{p-1}}l(r',R)\leq \frac{N-\ga p}{p-1}(\ln 2)^{-\frac{\eta}{p-1}}M^{\frac{1}{p-1}}$, thus in \eqref{L6} we can choose 
$$\ge=t^{-1}\norm{{\bf M}^\eta_{\ga p,R}[\gm_{B_{1}}]}_{L^\infty(\BBR^N)}^{\frac{1}{p-1}}=t^{-1}M^{\frac{1}{p-1}}~~\forall t >\max \{\ge_0^{-1},\frac{N-\ga p}{p-1}(\ln 2)^{-\frac{\eta}{p-1}} \}M^{\frac{1}{p-1}}$$
and as in the proof of \rprop{lambda},   $\left\{{\bf W}^R_{\ga,p}[\gm_{B_{1}}]>t\right\}\subset B_{2}$.\\
Then
\bel{B18}\BA{l}
\abs{\left\{{\bf W}^R_{\ga,p}[\gm_{B_1}]>3t\right\}\cap B_{2}}\leq c_0\exp\left(-\left(\frac{p-1-\eta}{4(p-1)}\right)^{\frac{p-1}{p-1-\eta}}\ga p\ln 2 M^{-\frac{1}{p-1-\eta}}t^{\frac{p-1}{p-1-\eta}}\right)\abs{B_{2}}.
\EA\ee
This can be written under the form 
\bel{B18}\BA{l}
\abs{\left\{F>t\right\}\cap B_{2}}\leq \abs{B_{2}}\chi_{(0,t_0]}+ c_0\exp\left(-\gd_0 t\right)\abs{B_{2}}\chi_{(t_0,\infty)}(t).
\EA\ee
where $F=M^{-\frac{1}{p-1-\eta}}\left({\bf W}^R_{\ga,p}[\gm_{B_1}]\right)^{\frac{p-1}{p-1-\eta}}$ and $t_0=\left(3\max \{\ge_0^{-1},\frac{N-\ga p}{p-1}(\ln 2)^{-\frac{\eta}{p-1}} \}\right)^{\frac{p-1}{p-1-\eta}}$.\\
Take $\gd\in (0,\gd_0)$, by Fubini's theorem
$$\BA {l}\myint{B_{2}}{}\exp\left(\gd F(x)\right) dx
=\gd \myint{0}{\infty}\exp\left(\gd t\right)
\abs{\{F>t\}\cap B_{2}}dt
\EA$$
Thus, $$\BA {l}
\myint{B_{2}}{}\exp\left(\gd F(x)\right) dx \leq \gd \myint{0}{t_0}\exp\left(\gd t\right) dt \abs{B_{2}}+c_0\gd 
\myint{t_0}{\infty}
\exp\left(-\left(\gd_0-\gd\right)t\right)dt\abs{B_{2}}
\\[4mm]\phantom{\myint{B_{2}}{}\exp\left(\gd F(x)\right) dx}
\leq \left(\exp\left(\gd t_0\right)-1\right) \abs{B_{2}}+\myfrac{c_0\gd}{\gd_0-\gd
} \abs{B_{2}}
\EA$$
which is the desired inequality.\qeda\\
\noindent{\Remark}By the proof of \rprop{lambda}, we see that $\ge_0 \geq \frac{c_{30}}{\max(1,\ln 40r )}$ where $c_{30}=c_{30}(N,\ga,p,\eta)>0$. Thus, $t_0\leq c_{31}\left(\max(1,\ln 40r )\right)^{\frac{p-1}{p-1-\eta}}$. Therefore $c_{29}\leq c_{32} \exp \left(c_{33} \left(\max(1,\ln 40r )\right)^{\frac{p-1}{p-1-\eta}}\right) $ where $c_{32}$ and $c_{33}$ depend on $N,\ga, p$ and $\eta$.
\subsection{Approximation of measures}
The next result is an extension of a classical result of Feyel and de la Pradelle \cite{FePra}. This type of result has been intensively used in the framework of Sobolev spaces since the pioneering work of Baras and Pierre \cite{BaPi}, but apparently it is new in the case of Bessel-Lorentz spaces. We recall that a sequence of bounded measures $\{\gm_n\}$ in $\Gw$ converges to some bounded measure $\gm$ in $\Gw$ in the {\it narrow topology} of $\mathfrak M^b(\Gw)$ if 
\bel{Na}
\lim_{n\to\infty}\myint{\Gw}{}\gf d\gm_n=\myint{\Gw}{}\gf d\gm\qquad\forall\gf\in C_b(\Gw):=C(\Gw)\cap L^{\infty}(\Gw).
\ee
\bth{upper} Assume $\Gw$ is an open subset of $\BBR^N$. Let $\ga>0$, $1<s<\infty$, $1\leq q<\infty$ and $\gm\in \mathfrak M_+(\Gw)$. If $\gm$ is absolutely continuous with respect to $C_{\ga,s,q}$ in $\Gw$, there exists a nondecreasing sequence
$\{\gm_n\}\subset \mathfrak M^b_+(\Gw)\cap (L^{\ga,s,q}(\BBR^N))'$, with compact support in $\Gw$ which converges to 
$\gm$ weakly in the sense of measures. Furthermore, if $\gm\in \mathfrak M^b_+(\Gw)$, then $\gm_n\rightharpoonup\gm$ in the narrow topology.
\es
\Proof {\it Step 1.} Assume that $\gm$ has compact support. Let $\gf\in L^{\ga,s,q}(\BBR^N)$ and $\tilde \gf$ its $C_{\ga,s,q}$-quasicontinuous representative. Since $\gm$ is abolutely continuous with respect to $C_{\ga,s,q}$, we can define the mapping
$$\gf\mapsto P(\gf)=\myint{\BBR^N}{}\tilde\gf^+d\gm\lfloor_{\Gw}
$$
where $\gm\lfloor_{\Gw}$ is the extension of $\gm$ by $0$ in $\Gw^c$.  By Fatou's lemma, $P$ is lower semicontinuous on $L^{\ga,s,q}(\BBR^N)$. Furthermore it is convex and potitively homogeneous of degree 1. If $Epi (P)$ denotes the epigraph of $P$, i.e.
$$Epi (P)=\{(\gf,t)\in L^{\ga,s,q}(\BBR^N)\ti \BBR:t\geq P(\gf)\},
$$
it is a closed convex cone.
Let $\ge>0$ and $\gf_0\in C^\infty_c$, $\gf_0\geq 0$. Since $(\gf_0,P(\gf_0)-\ge)\notin Epi (P)$, there exist $\ell\in (L^{\ga,s,q}(\BBR^N))'$, $a$ and $b$ in $\BBR$ such that 
\bel{C1}
a+bt+\ell(\gf)\leq 0\qquad\forall (\gf,t)\in Epi (P),
\ee
\bel{C2}a+b(P(\gf_0)-\ge)+\ell (\gf_0)>0.
\ee
Since $(0,0)\in Epi(P)$, $a\leq 0$. Since $(s\gf,st)\in Epi (P)$ for all $s>0$, $s^{-1}a+bt+\ell(\gf)\leq 0$, which implies 
$$bt+\ell(\gf)\leq 0\qquad\forall (\gf,t)\in Epi (P).$$
Finally, since $(0,1)\in Epi (P)$, $b\leq 0$. But if $b=0$ we would have $\ell(\gf)\leq -a$ for all $\gf\in L^{\ga,s,q}(\BBR^N)$. which would lead to $\ell=0$ and $a>0$ from (\ref{C2}), a contradiction. Therefore
$b<0$. Then, we put $\gth(\gf)=-\frac{\ell(\gf)}{b}$ and derive that, for any $(\gf,t)\in Epi (P)$, there holds  $\gth(\gf)\leq t$, and in particular
\bel{C3}
\gth(\gf)\leq P(\gf)\qquad\forall\gf\in L^{\ga,s,q}(\BBR^N). 
\ee
Since $\gf\leq 0\Longrightarrow P(\gf)=0$, $\gth$ is a positive linear functional on $L^{\ga,s,q}(\BBR^N)$. Furthermore
$$\sup_{\scriptsize\BA {l}\gf\in C^\infty_c(\BBR^N)\\
\norm\gf_{L^\infty}\leq 1\EA}\!\!\!\abs{\gth(\gf)}=\!\!\!\sup_{\scriptsize\BA {l}\gf\in C^\infty_c(\BBR^N)\\
\norm\gf_{L^\infty}\leq 1\EA}\!\!\!\gth(\gf)\leq \!\!\!\sup_{\scriptsize\BA {l}\gf\in C^\infty_c(\BBR^N)\\
\norm\gf_{L^\infty}\leq 1\EA}\!\!\!P(\gf)=P(1)=\gm(\Gw).
$$
By the Riesz representation theorem, there exists $\gs\in\mathfrak M_+(\BBR^N)$ such that
\bel{C4}\gth(\gf)=\myint{\BBR^N}{}\gf d\gs\qquad\forall\gf\in C_c^\infty(\BBR^N).
\ee
Inequality (\ref{C3}) implies $0\leq\gs\leq\gm\lfloor_{\Gw}$.  Thus $supp(\gs)\subset supp(\gm\lfloor_{\Gw})=supp(\gm)$ and $\gs$ vanishes on Borel subsets of  $C_{\ga,s,q}$ capacity zero, as $\gm$ does it, besides \eqref{C4} also values for all $\gf\in C^\infty(\BBR^N)$ . From (\ref{C2}), we have
$$\myint{\BBR^N}{}\tilde\gf_0 d\gs=\gth(\gf_0)>P(\gf_0)-\ge+\frac{a}{b}\geq \myint{\BBR^N}{}\tilde\gf_0d\gm\lfloor_{\Gw}-\ge.
$$
This implies 
\bel{C5}
0\leq \myint{\BBR^N}{}\tilde\gf_0 d(\gm\lfloor_{\Gw}-\gs)\leq\ge.
\ee
It remains to prove that $\gs\in (L^{\ga,s,q}(\BBR^N))'$. For all $f\in C^\infty_c(\BBR^N)$, $f\geq 0$, there holds
\bel{C6}
\myint{\BBR^N}{}{\bf G}_\ga [f] d\gs=\gth({\bf G}_\ga [f])\leq \norm{\gth}_{(L^{\ga,s,q}(\BBR^N))'}\norm{{\bf G}_\ga [f]}_{L^{\ga,s,q}(\BBR^N)},
\ee
since $\gth=-b^{-1}\ell$ and $\ell\in (L^{\ga,s,q}(\BBR^N))'$. Now, given $f\in L^{s,q}(\BBR^N)$, $f\geq 0$ and a sequence of molifiers $\{\gr_n\}$, $(\chi_{B_n}f)\ast\gr_n\in C^\infty_c(\BBR^N)$ and $(\chi_{B_n}f)\ast\gr_n\to f$ in $L^{s,q}(\BBR^N)$, where $\chi_{_{B_n}}$ is the indicator function of the ball $B_n$ centered at the origin of radius $n$. Furthermore, there is a subsequence 
$\{n_k\}$ such that $\lim_{n_k\to\infty}{\bf G}_\ga [(\chi_{B_{n_k}}f)\ast\gr_{n_k}](x)\to {\bf G}_\ga [f](x)$, $C_{\ga,s,q}$-quasi everywhere. Using Fatou's lemma and lower semicontinuity of the norm
$$\BA {l}
\myint{\BBR^N}{}{\bf G}_\ga [f] d\gs\leq\liminf_{n_k\to\infty}\myint{\BBR^N}{}{\bf G}_\ga [(\chi_{B_{n_k}}f)\ast\gr_{n_k}]d\gs\\[4mm]
\phantom{\myint{\BBR^N}{}{\bf G}_\ga [f] d\gs}
\leq \liminf_{n_k\to\infty}\norm{\gth}_{(L^{\ga,s,q}(\BBR^N))'}\norm{{\bf G}_\ga [(\chi_{B_{n_k}}f)\ast\gr_{n_k}]}_{L^{\ga,s,q}(\BBR^N)}\\[4mm]
\phantom{\myint{\BBR^N}{}{\bf G}_\ga [f] d\gs}
\leq \norm{\gth}_{(L^{\ga,s,q}(\BBR^N))'}\norm{{\bf G}_\ga [f]}_{L^{\ga,s,q}(\BBR^N)}.
\EA$$
Therefore (\ref{C6}) also holds  for all $f\in L^{s,q}(\BBR^N), f\geq 0$. 
Consequently $\gs\in \mathfrak M_+^b(\BBR^N)\cap (L^{\ga,s,q}(\BBR^N))'$ satisfies
\bel{C8}
\abs{\myint{\BBR^N}{}{\bf G}_\ga [f] d\gs}\leq \norm{\gth}_{(L^{\ga,s,q}(\BBR^N))'}\norm{{\bf G}_\ga [f]}_{L^{\ga,s,q}(\BBR^N)}\qquad\forall f\in L^{s,q}(\BBR^N).
\ee

\noindent {\it Step 2.} We assume that $\gm$ has no longer compact support. Set $\Gw_n=\{x\in \Gw:\dist(x,\Gw^c)\geq n^{-1},|x|\leq n\}$, then $\Gw_n\subset\overline{\Gw_n}\subset\Gw_{n+1}\subset\Gw$ for $n\geq n_0$ such that $\Gw_{n_0}\neq\emptyset$. Let $\{\gf_n\}\subset C^{\infty}_c(\BBR^N)$ be an increasing sequence such that $0\leq\gf_n\leq 1$, $\gf_n=1$ in a neighborhood of $\overline{\Gw_n}$ and $supp(\gf_n) \subset\Gw_{n+1}$. and let $\gn_n=\gf_n\gm$. For $n\geq n_0$ there is $\gs_n\in \mathfrak M_+^b(\BBR^N)\cap (L^{\ga,s,q}(\BBR^N))'$ with $0\leq\gs_n\leq\gn_n$ and
$$\frac{1}{n}>\myint{\Gw}{}\gf_n d(\gn_n-\gs_n)\geq \myint{\Gw_n}{} d(\gn_n-\gs_n)=\myint{\Gw_n}{} d(\gm-\gs_n).
$$
We set $\gm_n=\sup\{\gs_1,\gs_2,...,\gs_n\}$, then $\{\gm_n\}$ is nondecreasing and $supp(\gm_n)\subset\Gw_{n+1}$, and $\gm_n\in \mathfrak M_+^b(\BBR^N)\cap (L^{\ga,s,q}(\BBR^N))'$.
Finally, let $\gf\in C_c(\Gw)$ and $m\in\BBN^\ast$ such that $supp (\gf)\subset\Gw_m$. For all $n\geq m$, we have
$$\abs{\myint{\Gw}{}\gf d\gm_n-\myint{\Gw}{}\gf d\gm}\leq \abs{\myint{\Gw_n}{} d(\gm-\gm_n)}\norm\gf_{L^\infty(\BBR^N)}\leq\frac{1}{n}\norm\gf_{L^\infty(\BBR^N)}.
$$
Thus $\gm_n\rightharpoonup\gm$ weakly in the sense of measures.\medskip

\noindent {\it Step 3.} Assume that $\gm\in\mathfrak M^b_+(\Gw)$. Then $\gm_n(\Gw)\leq\gm(\Gw)$. Thus
$$\gm_n(\Gw)=\gm_n(\Gw_{n_0})+\displaystyle\sum_{k=n_0}^\infty\gm_n(\overline\Gw_{k+1}\setminus\Gw_{k})
$$
Since the sequence $\{\gm_n\}$ is nondecreasing and $\lim_{k\to\infty}\gm_n(\overline\Gw_{k+1}\setminus\Gw_{k})=\gm(\overline\Gw_{k+1}\setminus\Gw_{k})$by the previous construction, we 
obtain by monotone convergence
$$\lim_{n\to\infty}\gm_n(\Gw)=\gm(\Gw_{n_0})+\displaystyle\sum_{k=n_0}^\infty\gm(\overline\Gw_{k+1}\setminus\Gw_{k})=\gm(\Gw)
$$
Next we consider $\gf\in C_b(\Gw):=C(\Gw)\cap L^{\infty}(\Gw)$, then
$$\BA {l}
\abs{\myint{\Gw}{}\gf d\gm_n-\myint{\Gw}{}\gf d\gm}\leq \abs{\myint{\Gw}{}d(\gm-\gm_n)}
\norm\gf_{L^\infty(\Gw)}\leq (\gm(\Gw)-\gm_n(\Gw))\norm\gf_{L^\infty(\Gw)}\to 0.
\EA$$
Thus $\gm_n\rightharpoonup\gm$ in the narrow topology of measures.\qeda\medskip

As a consequence of \rth{upper} and \rth{equiv} we obtain the following.
\bth{upper+1} Let $p-1<s_1<\infty$,  $p-1<s_2\leq \infty$, $0<\ga p<N$, $R>0$ and $\gm\in \mathfrak M_+(\Gw)$. If $\gm$ is absolutely continuous with respect to the capacity 
$C_{\ga p,\frac{s_1}{s_1-p+1},\frac{s_2}{s_2-p+1}}$, there exists a nondecreasing sequence $\{\gm_n\}\subset \mathfrak M_+(\Gw)$ with compact support in $\Gw$ which converges to $\gm$ in the weak sense of measures and such that ${\bf W}_{\ga,p}^R[\gm_n]\in L^{s_1,s_2}(\BBR^N)$, for all $n$. Furthermore, if $\gm\in \mathfrak M^b_+(\Gw)$, $\gm_n$ converges to to $\gm$ in the narrow topology.
\es
\Proof By \rth{upper} there exists a nondecreasing sequence $\{\gm_n\}$ of nonnegative measures with compact support in $\Gw$, all elements of $ (L^{\ga p,\frac{s_1}{s_1-p+1},\frac{s_2}{s_2-p+1}}(\BBR^N))'$, which converges weakly to $\gm$. If $\gm\in \mathfrak M^b_+(\Gw)$, the convergence holds in the narrow topology. Noting that for a positive measure $\gs$ in $\BBR^N$,
$${\bf G}_{\ga p}[\gs]\in L^{\frac{s_1}{p-1},\frac{s_2}{p-1}}(\BBR^N)\Longleftrightarrow
\gs\in (L^{\ga p,\frac{s_1}{s_1-p+1},\frac{s_2}{s_2-p+1}}(\BBR^N))',
$$
it implies ${\bf G}_{\ga p}[\gm_n]\in L^{\frac{s_1}{p-1},\frac{s_2}{p-1}}(\BBR^N)$. Then, by \rth{equiv}, ${\bf W}_{\ga,p}^R[\gm_n]\in L^{s_1,s_2}(\BBR^N)$.\qeda
\section{Renormalized solutions}
\setcounter{equation}{0}
\subsection{Classical results}
Although the notion of renormalized solutions is becoming more and more present in the theory of quasilinear equations with measure data, it has not yet acquainted a popularity which could avoid us to present some of its main aspects. Let $\Gw$ be a bounded domain in $\BBR^N$.   If $\gm\in \mathfrak M^b(\Gw)$, we denote by $\gm^+$ and $\gm^-$ respectively its positive and negative part. We denote by $\mathfrak M_0(\Gw)$ the space of measures in $\Gw$ which are absolutely continuous with respect to the $c^{\Gw}_{1,p}$-capacity defined on a compact set $K\subset\Gw$ by
\bel{D1}
c^{\Gw}_{1,p}(K)=\inf\left\{\myint{\Gw}{}\abs{\nabla \gf}^pdx:\gf\geq \chi_K,\gf\in C^\infty_c(\Gw)\right\}.
\ee
We also denote $\mathfrak M_s(\Gw)$ the space of measures in $\Gw$ with support on a set of zero $c^{\Gw}_{1,p}$-capacity. Classically, any $\gm\in \mathfrak M^b(\Gw)$ can be written in a unique way under the form $\gm=\gm_0+\gm_s$ where $\gm_0\in \mathfrak M_0(\Gw)\cap\mathfrak M^b(\Gw)$ and $\gm_s\in \mathfrak M_s(\Gw)$.
We recall that any $\gm_0\in \mathfrak M_0(\Gw)\cap \mathfrak M^b(\Gw)$ can be written under the form $\gm_0=f-div \,g$ where $f\in L^1(\Gw)$ and $g\in L^{p'}(\Gw)$. \smallskip

For $k>0$ and $s\in\BBR$ we set $T_k(s)=\max\{\min\{s,k\},-k\}$. We recall that if $u$ is a measurable function defined and finite a.e. in $\Gw$, such that $T_k(u)\in W^{1,p}_0(\Gw)$ for any $k>0$, there exists a measurable function $v:\Gw\to \BBR^N$ such that $\nabla T_k(u)=\chi_{\abs u\leq k}v$ 
a.e. in $\Gw$ and for all $k>0$. We define the gradient $\nabla u$ of $u$ by $v=\nabla u$. We recall the definition of a renormalized solution given in \cite{DMOP}.

\bdef{renorm} Let $\gm=\gm_0+\gm_s\in\mathfrak M^b(\Gw)$. A measurable  function $u$ defined in $\Gw$ and finite a.e. is called a renormalized solution of 
\bel{D2}\BA {ll}
-\Gd _pu=\gm\qquad&\text{in }\Gw\\[0mm]\phantom{-\Gd _p}
u=0&\text{on }\prt\Gw,
\EA\ee
if $T_k(u)\in W^{1,p}_0(\Gw)$ for any $k>0$, $\abs{\nabla u}^{p-1}\in L^r(\Gw)$ for any $0<r<\frac{N}{N-1}$, and $u$ has the property that for any $k>0$ there exist $\gl_k^+,\gl_k^-\in \mathfrak M^b_+(\Gw)\cap \mathfrak M_0(\Gw)$, respectively concentrated on the sets $u=k$ and $u=-k$, with the property that 
$\gl_k^+\rightharpoonup\gm_s^+$, $\gl_k^-\rightharpoonup\gm_s^-$ in the narrow topology of measures, such that
\bel{D3}\BA {ll}
\myint{\{\abs u<k\}}{}\abs{\nabla u}^{p-2}\nabla u\nabla \gf\, dx=\myint{\{\abs u<k\}}{}\gf d\gm_0
+\myint{\Gw}{}\gf d\gl_k^+-\myint{\Gw}{}\gf d\gl_k^-,
\EA\ee
for every $\gf\in W^{1,p}_0(\Gw)\cap L^{\infty}(\Gw)$.
\es
\Remark If $u$ is a renormalized solution of problem (\ref{D2}) and $\gm\in\mathfrak M_+^b(\Gw)$, then $u\geq 0$ in $\Gw$. Indeed, taking $k>m>0$ and $\gf= T_m(\max\{-u,0\})$, then $0\leq \gf\leq m$ and we have
$$\BA {ll}
\myint{\{\abs u<k\}}{}\abs{\nabla u}^{p-2}\nabla u\nabla \gf dx =
\myint{\{\abs u<k\}}{}T_m(\max\{-u,0\}) d\gm_0+\myint{\Gw}{}T_m(\max\{-u,0\}) d\gl_k^+
\\[4mm]\phantom{\myint{\{\abs u<k\}}{}\abs{\nabla u}^{p-2}\nabla u\nabla \gf dx =
\myint{\{\abs u<k\}}{}T_m(\max\{-u,0\}) d\gm_0}
-\myint{\Gw}{}T_m(\max\{-u,0\}) d\gl_k^-
\\[4mm]\phantom{\myint{\{\abs u<k\}}{}\abs{\nabla u}^{p-2}\nabla u\nabla \gf dx}
\geq -m\gl_k^-(\Gw).
\EA$$
Thus
$$\myint{\Gw}{}\abs{\nabla T_m(\max\{-u,0\})}^p\leq m\gl_k^-(\Gw)
$$
Letting $k\to\infty$, we obtain $\nabla T_m(\max\{-u,0\})=0$ a.e., thus $u\geq 0$ 	a.e.  in $\Gw$. \medskip

We recall the following important results, see \cite[Th 4.1, Sec 5.1]{DMOP}. 
\bth{recall} 
\noindent  Let $\{\gm_n\}\subset \mathfrak M^b(\Gw)$ be a sequence such that $\sup_n\abs{\gm_n}(\Gw)<\infty$ and let $\{u_n\}$ be renormalized solutions of 
\bel{D4}\BA {ll}
-\Gd _pu_n=\gm_n\qquad&\text{in }\Gw\\[0mm]\phantom{-\Gd _p}
u_n=0&\text{on }\prt\Gw.
\EA\ee
Then, up to a subsequence, $\{u_n\}$ converges a.e. to a solution $u$ of $-\Gd _pu=\gm$ in the sense of distributions in $\Gw$, for some measure $\gm \in \mathfrak M^b(\Gw)$, and for every $k>0$, $k^{-1}\myint{\Gw}{}\abs{\nabla T_k(u)}^{p}\leq M$ for some $M>0$.
\es

Finally we recall the following fundamental stability result of \cite{DMOP} which extends \rth{recall}.

\bth{stab} Let $\gm=\gm_0+\gm_s^+-\gm_s^-\in \mathfrak M^b(\Gw)$, with $\gm_0=f-div\,g\in \mathfrak M_0(\Gw)$, $\gm_s^+,\gm_s^-\in \mathfrak M_s^+(\Gw)$. Assume there are sequences
$\{f_n\}\subset L^1(\Gw)$, $\{g_n\}\subset (L^{p'}(\Gw))^N$, $\{\eta_n^1\}, \{\eta_n^2\}\subset \mathfrak M_+^b(\Gw)$ such that $f_n\rightharpoonup f$ weakly in $L^1(\Gw)$, $g_n\to g$ in $L^{p'}(\Gw)$ and $div\,g_n$ is bounded in $\mathfrak M^b(\Gw)$, $\eta_n^1\rightharpoonup \gm_s^+$ and $\eta_n^2\rightharpoonup \gm_s^-$ in the narrow topology. If $\gm_n=f_n-div\,g_n+\eta_n^1-\eta_n^2$
and $u_n$ is a renormalized solution of (\ref{D4}), then, up to a subsequence, $u_n$ converges a.e. to a renormalized solution $u$ of (\ref{D2}). Furthermore $T_k(u_n)\to T_k(u)$ in $W^{1,p}_0(\Gw)$.
\es

\subsection{Applications}
We present below some interesting consequences of the above theorem.

\bcor{stab1} Let $\gm\in \mathfrak M^b(\Gw)$ with compact support in $\Gw$ and $\gw\in \mathfrak M^b(\Gw)$. Let $\{f_n\}\subset L^1(\Gw)$ which converges weakly to $f\in L^1(\Gw)$ and $\gm_n=\rho_n\ast\gm$ where $\{\rho_n\}$ is a sequence of mollifiers. If $u_n$ is a renormalized solution of 
\bel{D5}\BA {ll}
-\Gd _pu_n=f_n+\gm_n+\gw\qquad&\text{in }\Gw\\[0mm]\phantom{-\Gd _p}
u_n=0&\text{on }\prt\Gw,
\EA\ee
then, up to a subsequence, $u_n$ converges to a renormalized solution of 
\bel{D5}\BA {ll}
-\Gd _pu=f+\gm+\gw\qquad&\text{in }\Gw\\[0mm]\phantom{-\Gd _p}
u=0&\text{on }\prt\Gw.
\EA\ee
\es
\Proof We write $\gw=\tilde h-div\,\tilde g+\gw_s^+-\gw_s^-$ and $\gm= h-div\, g+\gm_s^+-\gm_s^-$, with $h,\tilde h\in L^1(\Gw)$, $g,\tilde g\in (L^{p'}(\Gw))^N$, $h$, $g$, $\gm_s^+$ and $\gm_s^-$ with support in a compact set $K\subset\Gw$. For $n_0$ large enough, $\rho_n\ast h$, $\rho_n\ast g$, $\rho_n\ast \gm_s^+$and $\rho_n\ast \gm_s^-$ have also their support in a fixed compact subset of $\Gw$ for all $n\geq n_0$. Moreover $\rho_n\ast h \to h$ and 
$\rho_n\ast g\to g$ in $L^1(\Gw)$ and $(L^{p'}(\Gw))^N$ respectively and $div\,\rho_n\ast g\to div\,g$ in $W^{-1,p'}(\Gw)$. Therefore
$$f_n+\gm_n+\gw=f_n+\tilde h+\rho_n \ast h-div\, (\tilde g+\rho_n\ast g)+\gw_s^++\rho_n\ast\gm_s^+-\gw_s^--
\rho_n\ast\gm_s^-$$
is an approximation of the measure $f+\gm+\gw$ in the sense of \rth{stab}. This implies the claim.\qeda

\bcor{stab2}Let $\gm_i\in \mathfrak M_+^b(\Gw)$, $i=1,2$, and $\{\gm_{i,n}\}\subset \mathfrak M_+^b(\Gw)$ be a nondecreasing and converging to $\gm_i$ in $\mathfrak M_+^b(\Gw)$. Let $\{f_n\}\subset L^1(\Gw)$ which converges to some $f$ weakly in $L^1(\Gw)$. Let $\{\gv _n\}\subset \mathfrak M^b(\Gw)$ which converges to some $\gv\in \mathfrak M_s(\Gw)$ in the narrow topology. For any $n\in\BBN$ let $u_n$ be a renormalized solution of 
\bel{D6}\BA {ll}
-\Gd _pu_n=f_n+\gm_{1,n}-\gm_{2,n}+\gv_n\qquad&\text{in }\Gw\\[0mm]\phantom{-\Gd _p}
u_n=0&\text{on }\prt\Gw.
\EA\ee
Then, up to a subsequence, $u_n$ converges a.e. to a renormalized solution of problem
\bel{D7}\BA {ll}
-\Gd _pu=f+\gm_{1}-\gm_{2}+\gv\qquad&\text{in }\Gw\\[0mm]\phantom{-\Gd _p}
u=0&\text{on }\prt\Gw.
\EA\ee
\es

The proof of  this results is based upon two lemmas

\blemma{lstab1} For any $\gm\in \mathfrak M_{0}(\Gw)\cap \mathfrak M_+^b(\Gw)$ there exists 
$f\in L^1(\Gw)$ and $h\in W^{-1,p'}(\Gw)$ such that $\gm=f+h$ and 
\bel{D8}
\norm{f}_{L^1(\Gw)}+\norm {h}_{W^{-1,p'}(\Gw)}+\norm{h}_{\mathfrak M^b(\Gw)}\leq 5\gm(\Gw).
\ee
\es
\Proof Following \cite {DaM} and the proof of \cite[Th 2.1]{BoGaOr}, one can write $\gm=\gf\gamma$ where 
$\gamma\in W^{-1,p'}(\Gw)\cap \mathfrak M_+^b(\Gw)$ and $0\leq\gf\in L^1(\Gw,\gamma)$. Let $\{\Gw_n\}_{n\in\BBN_\ast}$ be an increasing sequence of compact subsets of $\Gw$ such that $\cup_n \Gw_n=\Gw$. We define 
the sequence of measures $\{\gn_n\}_{n\in\BBN_\ast}$ by 
$$\BA {ll}
\gn_n=T_n(\chi_{\Gw_n}\gf)\gamma-T_{n-1}(\chi_{\Gw_{n-1}}\gf)\gamma\quad\text{for }n\geq 2\\
\gn_1=T_1(\chi_{\Gw_1}\gf)\gamma.
\EA$$ 
Since $\gn_k\geq 0$, then $\displaystyle\sum_{k=1}^\infty\gn_k=\gm$ with strong convergence  in $\mathfrak M^b(\Gw)$, $\norm{\gn_k}_{\mathfrak M^b(\Gw)}=\gn_k(\Gw)$ and 
$\displaystyle\sum_{k=1}^\infty\norm{\gn_k}_{\mathfrak M^b(\Gw)}=\gm(\Gw)$. Let $\{\rho_n\}$ be a sequence of mollifiers. We may assume that $\eta_n=\rho_n\ast \gn_n\in C^\infty_c(\Gw)$, 
$$\norm{\eta_n-\gn_n}_{W^{-1,p'}(\Gw)}\leq 2^{-n}\gm(\Gw)$$
Set $f_n=\displaystyle\sum_{k=1}^n\eta_k$, then 
$\displaystyle \norm{f_n}_{L^1(\Gw)}\leq \sum_{k=1}^n\norm{\eta_k}_{L^1(\Gw)}\leq \sum_{k=1}^n\norm{\gn_{k}}_{\mathfrak M^b(\Gw)}\leq \gm(\Gw).
$
If we define $f=\lim_{n\to\infty}f_n$, then $f\in L^1(\Gw)$ with $\norm{f}_{L^1(\Gw)}\leq \gm(\Gw)$. Set $h_n=\displaystyle\sum_{k=1}^n(\gn_k-\eta_k)$, then $h_n\in  W^{-1,p'}(\Gw)\cap \mathfrak M^b(\Gw)$, $\norm{h_n}_{W^{-1,p'}(\Gw)}\leq 2\gm(\Gw)$ and $h_n$ converges strongly in $W^{-1,p'}(\Gw)$ to some $h$ which satisfies $\norm{h}_{W^{-1,p'}(\Gw)}\leq 2\gm(\Gw)$. Since $\gm=f+h$ and $\norm{h}_{\mathfrak M^b(\Gw)}\leq 2\gm(\Gw)$, the result follows.\qeda

\blemma{lstab2} Let $\gm\in \mathfrak M_+^b(\Gw)$. If $\{\gm_n\}\subset \mathfrak M_+^b(\Gw)$ is a nondecreasing sequence which converges to $\gm$ in $\mathfrak M^b(\Gw)$, there exist $F_n,F\in L^1(\Gw)$, $G_n,G\in W^{-1,p'}(\Gw)$ and $\gm_{n\,s},\gm_s\in \mathfrak M_s(\Gw) $ such that
$$\gm_n=\gm_{n\,0}+\gm_{n\,s}=F_n+G_n+\gm_{n\,s}\quad\text{and }\;\gm=\gm_{0}+\gm_{s}=F+G+\gm_{s},
$$
such that $F_n\to F$ in $L^1(\Gw)$, $G_n\to G$ in $W^{-1,p'}(\Gw)$ and in $\mathfrak M^b(\Gw)$ and $\gm_{n\,s}\to\gm_{s}$ in $\mathfrak M^b(\Gw)$, and 
\es
\bel{D9}
\norm{F_n}_{L^1(\Gw)}+\norm {G_n}_{W^{-1,p'}(\Gw)}+\norm{G_n}_{\mathfrak M^b(\Gw)}+
\norm{\gm_{n\,s}}_{\mathfrak M^b(\Gw)}\leq 6\gm(\Gw).
\ee
\Proof Since $\{\gm_n\}$ is nondecreasing $\{\gm_{n\,0}\}$ and $\{\gm_{n\,s}\}$ share this property. Clearly 
$$\norm{\gm-\gm_n}_{\mathfrak M^b(\Gw)}=\norm{\gm_0-\gm_{n\,0}}_{\mathfrak M^b(\Gw)}+\norm{\gm_s-\gm_{n\,s}}_{\mathfrak M^b(\Gw)},$$
thus $\gm_{n\,0}\to \gm_{0}$ and $\gm_{n\,s}\to \gm_{s}$ in $\mathfrak M^b(\Gw)$. Furthermore
$\norm{\gm_{n\,s}}_{\mathfrak M^b(\Gw)}\leq\gm_s(\Gw)\leq\gm(\Gw)$. Set $\tilde\gm_{0\,0}=0$ and $\tilde\gm_{n\,0}=\gm_{n\,0}-\gm_{n-1\,0}$ for $n\in\BBN_\ast$. From \rlemma{lstab1}, for any $n\in\BBN$, one can find $f_n\in L^1(\Gw)$, $h_n\in W^{-1,p'}(\Gw)\cap \mathfrak M^b(\Gw)$ such that $\tilde\gm_{n\,0}=f_n+h_n$ and
$$\norm{f_n}_{L^1(\Gw)}+\norm{h_n}_{W^{-1,p'}(\Gw)}+\norm{h_n}_{\mathfrak M^b(\Gw)}\leq 5
\tilde\gm_{n\,0}(\Gw). 
$$
If we define $F_n=\displaystyle\sum_{k=1}^nf_k$ and $G_n=\displaystyle\sum_{k=1}^nh_k$, then 
$\gm_{n\,0}=F_n+G_n$ and 
$$\norm{F_n}_{L^1(\Gw)}+\norm{G_n}_{W^{-1,p'}(\Gw)}+\norm{G_n}_{\mathfrak M^b(\Gw)}\leq 5
\tilde\gm_{0}(\Gw). 
$$
Therefore the convergence statements and (\ref{D9}) hold.\qeda\medskip

\noindent{\it Proof of \rcor{stab2}}. We set $\gn_n=f_n+\gm_{n,1}-\gm_{n,2}+\gv_n$ and 
$\gn=f+\gm_{1}-\gm_{2}+\gv$. From \rlemma{lstab2} we can write 
$$\gn_n=f_n+F_{1\,n}-F_{2\,n}+G_{1\,n}-G_{2\,n}+\gm_{1\,n\,s}-\gm_{2\,n\,s}+\gv_n$$
and
$$\gn=f+F_{1}-F_{2}+G_{1}-G_{2}+\gm_{1\,s}-\gm_{2\,s}+\gv,
$$
and the convergence properties listed in the lemma hold. Therefore we can apply \rth{stab} and the conclusion follows.\qeda\medskip

 In the next result we prove the main pointwise estimates on renormalized solutions.
 \bth{pwest} Let $\Gw$ be a bounded domain of $\BBR^N$. Then there exists a constant $c>0$, dependent on $p$ and $N$ such that if $\gm\in \mathfrak M^b(\Gw)$ and $u$ is a renormalized solution of problem (\ref{D2}) there holds
 \bel{D10}
 -c{\bf W}^{2\,diam\,{\Gw}}_{1,p}[\gm^-]\leq u(x)\leq c{\bf W}^{2\,diam\,{\Gw}}_{1,p}[\gm^+]
 \quad\text{a.e. in }\Gw.\ee
 \es
 \Proof We claim the there exist renormalized solutions $u_1$ and $u_2$ of problem (\ref{D2}) with respective data  $\gm^+$ and $\gm^-$ such that
 \bel{D11}
 -u_2\leq u\leq u_1\qquad\text{a.e. in }\Gw.
\ee
We use the decomposition $\gm=\gm^+-\gm^-=(\gm_{0}^+-\gm_{s}^+)-(\gm_{0}^--\gm_{s}^-)$. We put 
$u_k=T_k(u)$, $\gm_k={\bf 1}_{\{\abs u<k\}}\gm_0+\gl_k^+-\gl_k^-$, $v_k={\bf 1}_{\{\abs u<k\}}\gm^+_0+\gl_k^+$.  Since $\gm_k\in \mathfrak M_0(\Gw)$, problem (\ref{D2}) with data $\gm_k$ admits a unique renormalized solution (see \cite{BoGaOr}), and clearly $u_k$ is such a solution. Since 
$v_k\in \mathfrak M_0(\Gw)$, problem (\ref{D2}) with data $v_k$ admits a unique solution $u_{k,1}$ which is furthermore nonnegative and dominates  $u_k$ a.e. in $\Gw$. From \rcor{stab2}, $\{u_{k,1}\}$ converges a.e. in $\Gw$ to a renormalized solution $u_1$ of (\ref{D2}) with data $\gm^+$ and $u\leq u_1$. Similarly $-u\leq u_2$ where $u_2$ is a renormalized solution of (\ref{D2}) with $\gm^-$. Finally, from \cite[Th 6.9]{PhVer} there is a positive constant $c$ dependent only on $p$ and $N$ such that 
 \bel{D12}
 u_1(x)\leq c{\bf W}^{2\,diam\,{\Gw}}_{1,p}[\gm^+]\quad\text{and }\;
  u_2(x)\leq c{\bf W}^{2\,diam\,{\Gw}}_{1,p}[\gm^-]
 \quad\text{a.e. in }\Gw.\ee
 This implies the claim.\qeda
 
 \section{Equations with absorption terms}
 
  \subsection{The general case}
  
  Let $g:\Gw\ti\BBR\mapsto \BBR$ be a Caratheodory function such that the map $s\mapsto g(x,s)$ is nondecreasing and odd for almost all $x\in\Gw$. If $U$ is a function defined in $\Gw$ we define the function $g\circ U$ in $\Gw$ by
  $$g\circ U(x)=g(x,U(x))\quad\text{for almost all }x\in\Gw.
  $$
  We consider the problem
   \bel{F1}\BA {ll}
-\Gd_pu+g\circ u=\gm\qquad&\text{in }\,\Gw\\\phantom{-\Gd_p+g\circ u}
u=0&\text{in }\,\prt\Gw.
\EA\ee
where $\gm\in \mathfrak M^b(\Gw)$. We say that $u$ is a {\it renormalized solution} of problem (\ref{F1}) if $g\circ u\in L^1(\Gw)$ and $u$ is a renormalized solution of 
   \bel{F2}\BA {ll}
-\Gd_pu=\gm-g\circ u\qquad&\text{in }\,\Gw\\\phantom{-\Gd_p}
u=0&\text{in }\,\prt\Gw.
\EA\ee
\bth{exist1} Let $\gm_i\in \mathfrak M_+^b(\Gw)$, $i=1,2$, such that there exists a nondecreasing sequences $\{\gm_{i,n}\}\subset \mathfrak M_+^b(\Gw)$, with compact support in $\Gw$,  converging to $\gm_i$ and $g\circ \left(c{\bf W}^{2\,diam\,{\Gw}}_{1,p}[\gm_{i,n}]\right)\in L^1(\Gw)$ with the same constant $c$ as in \rth{pwest}. Then there exists a renormalized solution of
   \bel{F3}\BA {ll}
-\Gd_pu+g\circ u=\gm_1-\gm_2\qquad&\text{in }\,\Gw\\\phantom{-\Gd_p+g\circ u}
u=0&\text{in }\,\prt\Gw,
\EA\ee
such that 
 \bel{F4}
 -c{\bf W}^{2\,diam\,{\Gw}}_{1,p}[\gm_2](x)\leq u(x)\leq c{\bf W}^{2\,diam\,{\Gw}}_{1,p}[\gm_1](x)
 \quad\text{a.e. in }\Gw.\ee
\es

\blemma{lexist1} Assume $g$ belongs to $L^\infty(\Gw\ti\BBR)$, besides the assumptions of \rth{exist1}. Let $\gl_i\in \mathfrak M_+^b(\Gw)$ ($i=1,2$), with compact support in $\Gw$. Then there exist renormalized solutions $u$, $u_i$, $v_i$ ($i=1,2$) to problems
   \bel{F5}\BA {ll}
-\Gd_pu+g\circ u=\gl_1-\gl_2\qquad&\text{in }\,\Gw\\\phantom{-\Gd_p+g\circ u}
u=0&\text{in }\,\prt\Gw,
\EA\ee
   \bel{F6}\BA {ll}
-\Gd_pu_i+g\circ u_i=\gl_i\qquad&\text{in }\,\Gw\\\phantom{-\Gd_p+g\circ u_i}
u_i=0&\text{in }\,\prt\Gw,
\EA\ee
   \bel{F7}\BA {ll}
-\Gd_pv_i=\gl_i\qquad&\text{in }\,\Gw\\\phantom{-\Gd_p}
v_i=0&\text{in }\,\prt\Gw,
\EA\ee
such that
   \bel{F8}\BA {ll}
-c{\bf W}^{2\,diam\,(\Gw)}_{1,p}[\gl_2](x)\leq -v_2(x)\leq -u_2(x)\leq u(x)\\[2mm]
\phantom{-c{\bf W}^{2\,diam\,(\Gw)}_{1,p}[\gl_2](x)----}
\leq u_1(x)\leq v_1(x)\leq c{\bf W}^{2\,diam\,(\Gw)}_{1,p}[\gl_1](x)
\EA\ee
for almost all $x\in\Gw$.
\es
\Proof Let $\{\gr_n\}$ be a sequence of mollifiers, $\gl_{i,n}=\gr_n\ast \gl_i,$ ($i=1,2$) and 
$\gl_n=\gl_{1,n}-\gl_{2,n}$. Then, for $n_0$ large enough,  $\gl_{1,n}$, $\gl_{2,n}$ and $\gl_n$ are bounded  with compact support in $\Gw$ for all $n\geq n_0$ and by minimization there exist unique solutions in $W^{1,p}_0(\Gw)$ to problems
$$\BA {ll}
-\Gd_pu_n+g\circ u_n=\gl_n\qquad&\text{in }\,\Gw\\\phantom{-\Gd_p+g\circ u_n}
u_n=0&\text{in }\,\prt\Gw,
\EA$$
$$\BA {ll}
-\Gd_pu_{i,n}+g\circ u_{i,n}=\gl_{i,n}\qquad&\text{in }\,\Gw\\\phantom{-\Gd_p+g\circ u_{i,n}}
u_{i,n}=0&\text{in }\,\prt\Gw,
\EA$$
$$\BA {ll}
-\Gd_pv_{i,n}=\gl_{i,n}\qquad&\text{in }\,\Gw\\\phantom{-\Gd_p}
v_{i,n}=0&\text{in }\,\prt\Gw,
\EA$$
and by the maximum principle, they satisfy
   \bel{F9}\BA {ll}
 -v_{2,n}(x)\leq -u_{2,n}(x)\leq u_n(x)\leq u_{1,n}(x)\leq v_{1,n}(x), \quad\forall x\in\Gw,\,\,\forall n\geq n_0.
\EA\ee
Since the $\gl_i$ are bounded measure and $g\in L^\infty(\Gw\ti\BBR)$ the the sequences of measures $\{\gl_{1,n}-\gl_{2,n}-g\circ u_n\}$, $\{\gl_{i,n}-g\circ u_{i,n}\}$ and  $\{\gl_{i,n}\}$ are uniformly bounded in $\mathfrak M^b(\Gw)$. Thus, by \rth{recall} there exists a subsequence, still denoted by the index $n$ such that $\{u_n\}$, $\{u_{i,n}\}$, $\{v_{i,n}\}$ converge a.e. in $\Gw$ to functions $\{u\}$, $\{u_{i}\}$, $\{v_i\}$ ($i=1,2$) when $n\to\infty$. Furthermore
$g\circ u_n$ and  $g\circ u_{i,n}$ converge in $L^1(\Gw)$ to $g\circ u$ and  $g\circ u_{i}$ respectively. By \rcor{stab1}, we can assume that $\{u\}$, $\{u_{i}\}$, $\{v_i\}$ are renormalized solutions of (\ref{F5})-(\ref{F7}), and by \rth {pwest}, $ v_{i}(x)\leq c{\bf W}^{2\,diam\,{\Gw}}_{1,p}[\gl_i](x)$, a.e. in $\Gw$. Thus we get (\ref{F8}).\qeda


\blemma{lexist2} Let $g$ satisfy the assumptions of \rth{exist1} and let $\gl_i\in \mathfrak M_+^b(\Gw)$ ($i=1,2$), with compact support in $\Gw$ such that $g\circ\left(c{\bf W}^{2\,diam\,(\Gw)}_{1,p}[\gl_i]\right)\in L^1(\Gw)$, where $c$ is the constant of \rth{exist1}. Then there exist renormalized solutions $u$, $u_i$ of the problems (\ref{F5})-(\ref{F6}) such that 
  \bel{F10}\BA {ll}
-c{\bf W}^{2\,diam\,(\Gw)}_{1,p}[\gl_2](x)\leq -u_2(x)\leq u(x)
\leq u_1(x)\leq c{\bf W}^{2\,diam\,(\Gw)}_{1,p}[\gl_1](x)
\EA\ee
for almost all $x\in\Gw$. Furthermore, if $\gw_i$, $\gth_i$ have the same properties as the $\gl_i$ and 
satisfy $\gw_i\leq\gl_i\leq\gth_i$, one can find solutions $u_{\gw_i}$ and $u_{\gth_i}$ of problems (\ref{F6}) with right-hand respective side $\gw_i$ and $\gth_i$, such that $u_{\gw_i}\leq u_i\leq u_{\gth_i}$.
\es
\Proof From \rlemma{lexist1} there exist renormalized solutions $u_n$, $u_{i,n}$ to problems
$$\BA {ll}
-\Gd_pu_n+T_n(g\circ u_n)=\gl_1-\gl_2\qquad&\text{in }\Gw\\\phantom{-\Gd_p+T_n(g\circ u_n)}
u_n=0&\text{on }\prt\Gw,
\EA$$
and
$$\BA {ll}
-\Gd_pu_{i,n}+T_n(g\circ u_{i,n})=\gl_i\qquad &\text{in }\Gw\\\phantom{-\Gd_p+T_n(g\circ u_{i,n})}
u_{i,n}=0&\text{on }\prt\Gw,
\EA$$
$i=1,2$, and they satisfy
  \bel{F11}\BA {ll}
-c{\bf W}^{2\,diam\,(\Gw)}_{1,p}[\gl_2](x)\leq -u_{2,n}(x)\leq u_n(x)
\leq u_{1,n}(x)\leq c{\bf W}^{2\,diam\,(\Gw)}_{1,p}[\gl_1](x).
\EA\ee
Since $\myint{\Gw}{}|g\circ u_{n}|dx\leq \gl_1(\Gw)+\gl_2(\Gw)$ and $\myint{\Gw}{}g\circ u_{i,n}dx\leq \gl_i(\Gw)$ thus as in \rlemma{lexist1} one can choose a subsequence, still denoted by the index $n$ such that 
$\{u_n, u_{1,n}, u_{2,n}\}$ converges a.e. in $\Gw$ to $\{u, u_{1}, u_{2}\}$ for which (\ref{F11}) is satisfied a.e. in $\Gw$. Since $g\circ\left(c{\bf W}^{2\,diam\,(\Gw)}_{1,p}[\gl_i]\right)\in L^1(\Gw)$ we derive from (\ref{F11}) and the dominated convergence theorem that $T_n(g\circ u_n)\to g\circ u$ and $T_n(g\circ u_{i,n})\to g\circ u_{i}$ in $L^1(\Gw)$. It follows from \rth{stab} that 
$u$ and $u_{i}$ are respective solutions of (\ref{F5}), (\ref{F6}). The last statement follows from the same assertion in \rlemma{lexist1}.\qeda\medskip

\noindent {\it Proof of \rth{exist1}.} From \rlemma {lexist2}, there exist renormalized solutions $u_n$, $u_{i,n}$ to problems
$$\BA {ll}
-\Gd_pu_n+g\circ u_n=\gm_{1,n}-\gm_{2,n}\qquad&\text{in }\Gw\\\phantom{-\Gd_p+g\circ u_n}
u_n=0&\text{on }\prt\Gw,
\EA$$
and
$$\BA {ll}
-\Gd_pu_{i,n}+g\circ u_{i,n}=\gm_{i,n}\qquad &\text{in }\Gw\\\phantom{-\Gd_p+g\circ u_{i,n}}
u_{i,n}=0&\text{on }\prt\Gw,
\EA$$
$i=1,2$ such that $\{u_{i,n}\}$ is nonnegative and nondecreasing and they satisfy
  \bel{F12}\BA {ll}
-c{\bf W}^{2\,diam\,(\Gw)}_{1,p}[\gm_2](x)\leq -u_{2,n}(x)\leq u_n(x)
\leq u_{1,n}(x)\leq c{\bf W}^{2\,diam\,(\Gw)}_{1,p}[\gm_1](x)
\EA\ee
a.e. in $\Gw$. As in the proof of \rlemma{lexist2}, up to the same subsequence, $\{u_{1,n}\}$,  $\{u_{2,n}\}$ and $\{u_n\}$ converge to $u_{1}$,  $u_{2}$ and $u$ a.e. in $\Gw$. Since $g\circ u_{i,n}$ are nondecreasing, positive and $\myint{\Gw}{}g\circ u_{i,n}dx\leq \gm_{i,n}(\Gw)\leq \gm_{i}(\Gw)$, it follows from the monotone convergence theorem that $\{g\circ u_{i,n}\}$ converges to $g\circ u_{i}$ in $L^1(\Gw)$. Finally, since $\abs{g\circ u_{n}}\leq g\circ u_{1}+g\circ u_{2}$, $\{g\circ u_{n}\}$ converges to $g\circ u$ in $L^1(\Gw)$ by dominated convergence. Applying \rcor{stab2} we conclude that $u$ is a renormalized solution of (\ref{F3}) and that (\ref{F4}) holds.\qeda\medskip

\subsection{Proofs of \rth{power} and \rth{exp}}

We are now in situation of proving the two theorems stated in the introduction.\medskip

\noindent{\it Proof of \rth{power}.} {\bf 1}- Since $\gm$ is absolutely continuous with respect to the capacity 
$C_{p,\frac{Nq}{Nq-(p-1)(N-\gb))},\frac{q}{q+1-p}}$, $\gm^+$ and $\gm^-$ share this property. By \rth{upper+1} there exist two nondecreasing sequences $\{\gm_{1,n}\}$ and $\{\gm_{2,n}\}$ of positive bounded measures with compact support in $\Gw$ which converge to $\gm^+$ and $\gm^-$ respectively and which have the property that ${\bf W}^{R}_{1,p}[\gm_{i,n}]\in L^{\frac{Nq}{N-\gb},q}(\BBR^N)$, for $i=1,2$ and all $n\in\BBN$. Furthermore, with $R=diam\,(\Gw)$,
  \bel{F13}\BA {ll}
\myint{\BBR^N}{}\frac{1}{\abs x^\gb}\left({\bf W}^{2R}_{1,p}[\gm_{i,n}](x)\right)^qdx
\leq \myint{0}{\infty}\left(\frac{1}{\abs {.}^\gb}\right)^\ast(t)\left(\left({\bf W}^{2R}_{1,p}[\gm_{i,n}]\right)^\ast(t)\right)^qdt
\\[4mm]\phantom{\myint{\BBR^N}{}\frac{1}{\abs x^\gb}\left({\bf W}^{2R}_{1,p}[\gm_{i,n}](x)\right)^qdx}
\leq c_{34}\myint{0}{\infty}\myfrac{1}{t^{\frac{\gb}{N}}}\left(\left({\bf W}^{2R}_{1,p}[\gm_{i,n}]\right)^\ast(t)\right)^qdt
\\[4mm]\phantom{\myint{\BBR^N}{}\frac{1}{\abs x^\gb}\left({\bf W}^{2R}_{1,p}[\gm_{i,n}](x)\right)^qdx}
\leq c_{34}\norm{{\bf W}^{2R}_{1,p}[\gm_{i,n}]}_{L^{\frac{Nq}{N-\gb},q}(\BBR^N)}^q\\\phantom{\myint{\BBR^N}{}\frac{1}{\abs x^\gb}\left({\bf W}^{2R}_{1,p}[\gm_{i,n}](x)\right)^qdx}
<\infty.
\EA\ee
Then the result follows from \rth{exist1}.\smallskip

\noindent {\bf 2}- Because $\gm$ is absolutely continuous with respect to the capacity 
$C_{p,\frac{Nq}{Nq-(p-1)(N-\gb))},1}$, so are $\gm^+$ and $\gm^-$. Applying again \rth{upper+1} there exist two nondecreasing sequences $\{\gm_{1,n}\}$ and $\{\gm_{2,n}\}$ of positive bounded measures with compact support in $\Gw$ which converge to $\gm^+$ and $\gm^-$ respectively and such that ${\bf W}^{R}_{1,p}[\gm_{i,n}]\in L^{\frac{Nq}{N-\gb},1}(\BBR^N)$. This implies in particular 
\bel{F14}
\left({\bf W}^{2R}_{1,p}[\gm_{i,n}](.)\right)^\ast(t)\leq c_{35}t^{-\frac{N-\gb}{Nq}},\qquad\forall t>0,
\ee
for some $c_{34}>0$. Therefore, by \rth{equiv}
  \bel{F15}\BA {ll}
\myint{\Gw}{}\frac{1}{\abs x^\gb}g\left(c{\bf W}^{2R}_{1,p}[\gm_{i,n}](x)\right)dx
\leq \myint{0}{\abs\Gw}\left(\frac{1}{\abs {.}^\gb}\right)^\ast(t)g\left(c\left({\bf W}^{2R}_{1,p}[\gm_{i,n}]\right)^\ast(t)\right)dt
\\[4mm]\phantom{\myint{\Gw}{}\frac{1}{\abs x^\gb}g\left(c{\bf W}^{2R}_{1,p}[\gm_{i,n}](x)\right)dx}
\leq c_{36}\myint{0}{\abs\Gw}\myfrac{1}{t^{\frac{\gb}{N}}}g\left(c\left({\bf W}^{2R}_{1,p}[\gm_{i,n}]\right)^\ast(t)\right)dt
\\[4mm]\phantom{\myint{\BBR^N}{}\frac{1}{\abs x^\gb}g\left(c{\bf W}^{2R}_{1,p}[\gm_{i,n}](x)\right)dx}
\leq c_{36}\myint{0}{\abs\Gw}\myfrac{1}{t^{\frac{\gb}{N}}}g\left(c_{35}ct^{-\frac{N-\gb}{Nq}}\right)dt
\\\phantom{\myint{\BBR^N}{}\frac{1}{\abs x^\gb}g\left(c{\bf W}^{2R}_{1,p}[\gm_{i,n}](x)\right)dx}
\leq c_{37} \myint{a}{\infty}g(t)t^{-q-1}dt
\\\phantom{\myint{\BBR^N}{}\frac{1}{\abs x^\gb}g\left(c_{22}{\bf W}^{2R}_{1,p}[\gm_{i,n}](x)\right)dx}
<\infty,
\EA\ee
where $a>0$ depends on $\abs\Gw$, $c_{35}c$, $N$, $\gb$, $q$.  Thus the result follows by \rth{exist1}.
\qeda\medskip

\noindent{\it Proof of \rth{exp}.} Again we take $R=diam\,(\Gw)$. Let $\{\Gw_n\}_{n\in\BBN_\ast}$ be an increasing sequence of compact subsets of $\Gw$ such that $\cup_n \Gw_n=\Gw$. We define $\gm_{i,n}=T_n(\chi_{\Gw_n}f_i)+\chi_{\Gw_n}\gn_i$ ($i=1,2$). Then $\{\gm_{1,n}\}$ and $\{\gm_{2,n}\}$ are nondecreasing sequences of elements of $\mathfrak M_+^b(\Gw)$ with compact support, and they converge to $\gm^+$ and $\gm^-$ respectively. Since for any $\ge>0$ there exists $c_\ge>0$ such that 
\bel{F16}
\left({\bf W}^{2R}_{1,p}[\gm_{i,n}]\right)^\gl
\leq c_\ge n^{\frac{\gl}{p-1}}+(1+\ge)\left({\bf W}^{2R}_{1,p}[\gn_i]\right)^\gl,
\ee
a.e. in $\Gw$, it follows
\bel{F17}
\exp\left(\gt\left(c{\bf W}^{2R}_{1,p}[\gm_{i,n}]\right)^\gl\right)
\leq c_{\ge,n,c} \exp\left(\gt(1+\ge)\left(c{\bf W}^{2R}_{1,p}[\gn_i]\right)^\gl\right).
\ee
If there holds 
\bel{F18}
\norm{{\bf M}^{\frac{(p-1)(\gl-1)}{\gl}}_{p,2R}[\gn_i]}_{L^\infty(\Gw)}<
\left(\myfrac{p\ln 2}{\gt (12\gl c)^\gl}\right)^{\frac{p-1}{\gl}},
\ee
we can choose $\ge>0$ small enough so that
$$\gt (1+\ge)c^\gl<\myfrac{ p\ln 2}{(12\gl)^{\gl}\norm{{\bf M}^{\frac{(p-1)(\gl-1)}{\gl}}_{p,2R}[\gn_i]}^{\frac{\gl}{p-1}}_{L^\infty(\Gw)}}.
$$
Hence, by \rth{exp1} with $\eta=\frac{(p-1)(\gl-1)}{\gl}$, 
$\exp\left(\gt(1+\ge)\left(c{\bf W}^{2R}_{1,p}[\gn_i]\right)^\gl\right)\in L^1(\Gw)$, which implies
$\exp\left(\gt\left(c{\bf W}^{2diam\,(\Gw)}_{1,p}[\gm_{i,n}]\right)^\gl\right)\in L^1(\Gw).$
We conclude by \rth{exist1}.
\qeda


\begin{thebibliography}{99}

\bibitem{AdHe} D. R. Adams and L. I. Hedberg: {\em Function Spaces and Potential Theory}, Springer, New York, (1996).

\bibitem{ArMuSz} N. Aronszjan, P. Mulla, P. Szeptycki: {\em On spaces of potentials connected with $L^q$ classes}, {\bf Ann. Inst. Fourier Grenoble 13}, 211-306 (1963).

\bibitem{BaPi} P. Baras, M. Pierre: {\em Singularit\'es \'eliminables pour des \'equations semi lin\'eaires},
 {\bf Ann. Inst. Fourier Grenoble 34}, 185-206 (1984).
\bibitem{BeBr} Ph. Benilan, H. Brezis: Nonlinear problems related to the Thomas-Fermi equation, {\bf unpublished paper}, see \cite{Br}
\bibitem {Bi1}M. F. Bidaut-Veron: {\em  Necessary conditions of existence for an
elliptic equation with source term and measure data involving p-Laplacian},
Proc. 2001 Luminy Conf. on Quasilinear Elliptic and Parabolic Equations and
Systems, {\bf Elect. J. Diff. Equ. Conf. 8}, 23-34 (2002).

\bibitem {Bi2}M. F. Bidaut-Veron: {\em Removable singularities and existence for a quasilinear equation with absorption or source term and measure data}, {\bf  Adv. Nonlinear Stud. 3}, 25Ð63 (2003).

\bibitem {BoGaOr}L. Boccardo, T. Galouet, L. Orsina: {\em Existence and uniqueness of entropy solutions for nonlinear elliptic equations with measure data}, {\bf  Ann. Inst. H. Poincar\'e, Anal. Non Lin\'eaire 13}, 539-555 (1996).

\bibitem{Br} H. Brezis, Some variational problems of the Thomas-Fermi type, in Variational Inequalities, eds. R.W. Cottle, F. Giannessi and J. L. Lions, Wiley, Chichester (1980), 53-73.

\bibitem{DaM} G. Dal Maso: {\em On the integral representation of certain local functionals}, {\bf Ricerche Mat. 32}, 85-113 (1983).

\bibitem {DMOP} G. Dal Maso, F. Murat, L. Orsina, A. Prignet: {\em Renormalized solutions of elliptic equations with general measure data}, {\bf  Ann.Sc. Norm. Sup. Pisa 28}, 741-808 (1999)

\bibitem{FePra} D. Feyel, A. de la Pradelle: {\em Topologies fines et compactifications associ\'ees \`a certains espaces de Dirichlet}, {\bf Ann. Inst. Fourier Grenoble 27}, 121-146 (1977).


\bibitem{Gr} L. Grafakos: {\em Classical Fourier Analysis} 2nd ed., Graduate Texts in Math. {\bf 249}, Springer-Verlag (2008).


\bibitem{HeKiMa} J. Heinonen, T. Kilpelainen, O. Martio {\em Nonlinear Potential Theory}, Oxford Univ. Press, Oxford (1993).

\bibitem{HoJa} P. Honzik and B. Jaye: {\em On the good-$\gl$ inequality for nonlinear potentials}, {\bf Proc. Amer. Math. Soc. 140}, 4167-4180 (2012).

\bibitem{MuWh} B. Muckenhoupt, R. Wheeden: {\em Weighted norm inequalities for fractional integrals}, {\bf Trans. Amer. Math. Soc. 192}, 261-274 (1974).

\bibitem{Nei} R. O. Neil, {\em Convolution operators on $L^{p,q}$ spaces}, {\bf Duke Math. J. 30}, 129-142 (1963).

\bibitem{PhVer} N. C. Phuc, I. E. Verbitsky: {\em Quasilinear and Hessian equations of Lane-Emden type},{\bf Ann. Math. 168}, 859-914 (2008).

\bibitem{St} E. M. Stein: {\em Singular Integrals and Differentiability of Functions}, Princeton Univ. Press, Princeton N.J. (1971).

\bibitem{Tur} B. O. Tureson, {\em Nonlinear Potential Theory and Sobolev Spaces}, Springer-Verlag (2000).
%
\bibitem{V1} L. V\'eron: {\em Elliptic equations involving measures}, Stationary partial differential equations. Vol. I, 593--712, {\bf Handb. Differ. Equ.}, North-Holland, Amsterdam, (2004).
%
\bibitem{V2} L.V\'eron: {\em Singularities of solutions of second other Quasilinear Equations}, Pitman Research Notes in Math. Series 353, Adison Wesley, Longman 1996.

\bibitem{Zie} W. Ziemer: {\em Weakly Differentiable Functions}, Springer-Verlag (1989).
\end{thebibliography}
\end{document}